\documentclass[10pt,a4paper]{amsart}
\usepackage[latin1]{inputenc}
\usepackage[english]{babel}
\usepackage{amsmath}
\usepackage{amsfonts}
\usepackage{amssymb}
\usepackage{epsfig}
\usepackage[left=2cm,right=2cm,top=2cm,bottom=2cm]{geometry}

\newcommand{\C}{\mathbb{C}}
\newcommand{\R}{\mathbb{R}}
\newcommand{\K}{\mathbb{K}}
\newcommand{\Q}{\mathbb{Q}}
\renewcommand{\H}{\mathcal{H}}
\newcommand{\N}{\mathbb{N}}
\newcommand{\z}{\zeta}
\renewcommand{\sp}[2]{\big( {#1} \mid {#2} \big)}
\newcommand{\dom}{D}
\newcommand{\norm}[1]{\left\Vert{#1}\right\Vert}
\newcommand{\abs}[1]{\left|{#1}\right|}

\newcommand{\B}{\mathcal{B}}

\renewcommand{\Re}{\operatorname{Re}\,}

\newcommand{\He}{\operatorname{Sym}\,}
\newcommand{\Sym}{\He}
\newcommand{\Lip}{\operatorname{Lip}}
\newcommand{\dd}{\, \mathrm{d}}
\newcommand{\ee}{\mathrm{e}}
\newcommand{\ii}{\, \mathrm{i}}

\newtheorem{thrm}{Theorem}[section]
\newtheorem{dfntn}[thrm]{Definition}
\newtheorem{rmrk}[thrm]{Remark}

\newtheorem{prpstn}[thrm]{Proposition}
\newtheorem{assumption}[thrm]{Assumption}

\newtheorem{lmm}[thrm]{Lemma}

\begin{document}

\author{Bj\"orn Augner}
\address{Technische Universit\"at Darmstadt, Fachbereich Mathematik, Mathematische Modellierung und Analysis, Schlossgartenstr.\ 7, 64289 Darmstadt.}
\email{augner@mma.tu-darmstadt.de}
\thanks{This work has been partly supported by Deutsche Forschungsgemeinschaft (Grant JA 735/8-1).}
\title[Exponential Stabilisation of Serially Connected Euler-Bernoulli Beams]{Uniform Exponential Stabilisation of Serially Connected Inhomogeneous Euler-Bernoulli Beams}

\subjclass[2010]{93D15, 35B35, 35G46, 34B09}
\keywords{Euler-Bernoulli beam, inhomogeneous distributed parameter systems, serially connected beams, exponential stabilisation, frequency domain method.} 

 \begin{abstract}
  We consider a chain of Euler-Bernoulli beams with spatial dependent mass density, modulus of elasticity and area moment which are interconnected in dissipative or conservative ways and prove uniform exponential energy decay of the coupled system for suitable dissipative boundary conditions at one end and suitable conservative boundary conditions at the other end.
  We thereby generalise some results of G.~Chen, M.C.~Delfour, A.M.~Krall and G.~Payre from the 1980's to the case of spatial dependence of the parameters.
 \end{abstract}
 

 \maketitle
 
 \section{Introduction}
 Beam equations became the focus of attention for mathematical modelling, analysis and numerics of complex, multi-component systems, in particular mechanical systems for the modelling of aeroplanes, bridges, nowadays more and more electromechanical systems and robotics, at least since the 1980's.
 Several types of partial differential equations serve as, and compete as models for such vibrating beams or strings: from the wave equation, probably one of the most commonly and most detailed discussed models in mathematics, to the Rayleigh beam and the Euler-Bernoulli beam, to the Timoshenko beam and even more sophisticated models. In many cases, these equations are non-linear in principle, but for the analysis and numerics of complex systems it is often useful, to consider the linear or linearised versions of these equations.
 In this article, we treat the linear Euler-Bernoulli beam model
  \[
   \rho(\z) \omega_{tt}(t,\z) 
    + (EI(\z) \omega_{\z\z})_{\z\z}(t,\z)
    = 0,
    \quad
    t \geq 0, \, \z \in (0,l)
  \]
 where $\rho(\z)$ denotes the \emph{mass density times cross section area} of a beam of length $l > 0$, and $E(\z)$ and $I(\z)$ its \emph{modulus of elasticity} and \emph{area moment of the cross section}, respectively.
 G.~Chen and several coauthors \cite{ChenEtAl_1987}, \cite{ChenEtAl_1987a}, \cite{ChenEtAl_1989} considered three particular important situations for the Euler-Bernoulli beam:
  \begin{enumerate}
   \item
    A single beam is stabilised by dissipative boundary feedback at one end of the beam and conservative boundary conditions at the other end \cite{ChenEtAl_1987a}.
   \item
    A pair of \emph{identical} beams is damped via dissipative point feedback at the joint \cite{ChenEtAl_1989}.
   \item
    An arbitrary long, but finite chain of serially connected beams is damped at one end of the chain \cite{ChenEtAl_1987}.
  \end{enumerate}
 In all these cases, the authors assume that the beam parameters $\rho$, $E$ and $I$ are constant along each of the beams.
 Since then for all three cases the corresponding articles inspired further mathematical research for more general models.
 E.g.\ in \cite{GuoHuang_2004} and \cite{Morgul_2001}, situations have been considered where for a single Euler-Bernoulli beam the collocated feedback at the dissipative end is perturbed, i.e.\ the feedback input cannot be expressed solely by (the traces of) the energy variables $\omega_{t}$ and $EI \omega_{\z\z}$ and their spatial derivatives.
 \newline
 Other works, e.g.\ \cite{AmmariTucsnak_2000}, \cite{GuoChan_2001}, \cite{AmmariLiuTucsnak_2002}, \cite{GuoXie_2004}, \cite{AbriolaEtAl_2017} further dealt with the problem of dissipative point feedback at the joint between two Euler-Bernoulli beams.
 These works highlighted that such a feedback law is not a good choice for exponential stabilisation (or, it is not a good model for such systems), because typically they gave the result that the property of asymptotic and uniform exponential stability depends on whether the actuation position $\xi \in (0,l)$, at which the damper acts, the fraction $\frac{\xi}{l} \in (0,1)$ lies in some subset $\tilde Q$ of $\Q \cap [0,1]$ which is still dense in $[0,1]$.
 A very unsatisfactory result from engineering perspective.
 \newline
 At the same time, more general networks defining the interconnection structure of Euler-Bernoulli beams and their stability properties have been considered, e.g.\ in \cite{DekoninckNicaise_2000}, \cite{MercierRegnier_2008} and \cite{MercierRegnier_2008a}.
 \newline
 The methods used mostly for proving stability essentially break down into three more or less heavily used methods:
  \begin{enumerate}
   \item
    Construction of a suitable \emph{Lyapunov function}:
    This method has been applied in \cite{ChenEtAl_1987}, \cite{AmmariTucsnak_2000}.
   \item
    Analysis of the asymptotic behaviour of the (discrete) eigenvalues $\lambda_n$ for $n \rightarrow \infty$, see e.g.\ \cite{ChenEtAl_1987a}, \cite{DekoninckNicaise_2000}, \cite{GuoChan_2001}, \cite{GuoXie_2004}, \cite{MercierRegnier_2008}, \cite{MercierRegnier_2008a}.
   \item
    \emph{Frequency domain method}: Resolvent estimates on the imaginary axis based on the Gearhart-Pr\"uss-Huang Theorem, i.e.\
     \[
      \sup_{\beta \in \R} \norm{(\ii \beta - \mathcal{A})^{-1}} < \infty,
     \]
    e.g.\ in \cite{Morgul_2001}, \cite{AmmariLiuTucsnak_2002}, \cite{GuoHuang_2004}.
  \end{enumerate}
 Each of these methods has its own advantages and disadvantages.
 E.g.\ the first method is suitable to allow for non-linear perturbations in the dissipative boundary feedback, but the method seems to be restricted to Euler-Bernoulli beams with almost homogeneous parameters $\rho$ and $EI$, cf.\ \cite{Augner_2018+}, and it is not clear at all whether all cases for which uniform exponential stabilisation is already known can be covered by this method as well.
 \newline
 Even more restrictive seems the second method, which mainly can be used for homogeneous beam models, whereas the frequency domain method in generally is suitable for non-homogeneous beams as well (and will be applied in this article).
 At the same time, both the second and the third method are restricted to the case of linear boundary feedback, and leave stability questions concerning nonlinear feedback wide open.
 \newline
 Note that the papers listed above almost exclusively cover \emph{homogeneous} beam equations, i.e.\ $\rho$ and $EI$ are constant, at least on each beam.
 This brings up the question:
 Is homogeneity of the beams only a technical restriction for the proofs? Can the general inhomogeneous case be reduced to the special homogeneous case? Does a (sufficiently regular) inhomogeneity influence well-posedness or stability at all?
 As it turns out, for the last question, which actually consists of two separate questions (well-posedness and stability), one of which has an easy answer, the other not so.
In fact, for dissipative systems well-posedness (in the sense of semigroup generation, i.e.\ existence, uniqueness and continuous dependence on the initial datum for the corresponding abstract Cauchy problem) is invariant under perturbation by a coercive and continuous operator, see e.g.\ \cite[Lemma 7.2.3]{JacobZwart_2012} or the much more general results in \cite{CalvertGustafson_1972}.
  (For a background on strongly continuous semigroups ($C_0$-semigroups), we refer to the monograph \cite{EngelNagel_2000}.)
  Does the same result hold if the term \emph{strongly continuous contraction semigroup} is replaced by \emph{uniformly exponentially stable, strongly continuous contraction semigroup}?
  Unfortunately not! Actually, there are already examples on finite dimensional Hilbert spaces which serve as counter examples, and in the class of \emph{infinite-dimensional port-Hamiltonian systems} \cite{VanDerSchaftMaschke_2002}, \cite{LeGorrecZwartMaschke_2005}, \cite{JacobZwart_2012}, i.e.\ a hyperbolic vector-valued PDE on an interval (in which form the Euler-Bernoulli beam can be rewritten) a striking counter example is known \cite{Engel_2013}.
  Though this particular counter example does not belong to the class of Euler-Bernoulli beams, yet it motivates the standpoint we take in this paper: Stability of inhomogeneous beams should be addressed additionally to the question of stability for their homogeneous counterparts. Therefore, we generalise the results of \cite{ChenEtAl_1987} in this direction, which -- to our knowledge -- has not yet been achieved up to now.
G.~Chen et al.\ \cite{ChenEtAl_1987} investigated a system of Euler-Bernoulli beams which are serially interconnected (in a conservative or dissipative way), and which is damped at one of the two ends of the chain, e.g.\
{\allowdisplaybreaks[1]
 \begin{align*}
  \rho^j \omega_{tt}(t,\z)
   + (E^j I^j \omega_{\z\z})_{\z\z}(t,\z)
   &= 0,
   &&t \geq 0, \, \z \in (l^{j-1}, l^j), \, j = 1, \ldots, m
   \\
  \omega(t,0)
   &= 0,
   &&t \geq 0
   \\
  \omega_\z(t,0)
   &= 0,
   &&t \geq 0
   \\
  \omega(t, l^j -)
   &= \omega(t, l^j +),
   &&t \geq 0, \,j = 1, \ldots, m-1
   \\
  \omega_\z(t,l^j -)
   &= \omega_\z(t, l^j +),
   &&t \geq 0, \, j = 1, \ldots, m-1
   \\ 
  - (E^j I^j \omega_{\z\z}) (t, l^j -)
   &= - (E^{j+1} I^{j+1} \omega_{\z\z}) (t, l^j +),
   &&t \geq 0, \, j = 1, \ldots, m-1
   \\
  (E^j I^j \omega_{\z\z})_\z (t,l^j -)
   &= (E^j I^j \omega_{\z\z})_\z (t, l^j +),
   &&\, j = 1, \ldots, m-1
   \\
  - (E^m I^m \omega_{\z\z})(t,L)
   &=0,
   &&t \geq 0
   \\
  (E^m I^m \omega_{\z\z})_\z(t,L)
   &= \kappa \omega_t(t,L)
   &&t \geq 0
   \\
  \omega(0,\z)
   &= \omega_0(\z),
   &&\z \in (l^{j-1}, l^j), \, j = 1, \ldots, m
   \\
  \omega_t(0,\z)
   &= \omega_1(\z),
   &&\z \in (l^{j-1}, l^j), \, j = 1, \ldots, m
 \end{align*}
}where $0 = l^0 < l^1 < \ldots < l^m = L$ is a division of the interval $(0, L)$ for some $L > 0$ and $m \in \N$, and $\kappa > 0$ is some damping parameter.
Here and in the following, we write
 \[
  f(\zeta \pm)
   := \lim_{\omega \rightarrow \zeta \pm} f(\omega)
 \]
for the one-sided limits of a function at position $\z$.
For the special case of a pair of beams ($m = 2$) of unit total length ($L = 1$) and a joint at position $l^1 = l \in (0,1)$, this system reads as
{\allowdisplaybreaks[1]
 \begin{align*}
  \rho^1 \omega_{tt}(t,\z) + (E^1 I^1 \omega_{\z\z})_{\z\z}(t,\z)
   &= 0,
   &&t \geq 0, \, \z \in (0,l)
   \\
  \rho^2 \omega_{tt}(t,\z) + (E^2 I^2 \omega_{\z\z})_{\z\z}(t,\z)
   &= 0,
   &&t \geq 0, \, \z \in (l,1)
   \\
  \omega(t,0)
   &= 0,
   &&t \geq 0
   \\
  \omega_\z(t,0)
   &= 0,
   &&t \geq 0
   \\
  \omega_\z(t,l-)
   &= \omega_\z(t,l+),
   &&t \geq 0
   \\
  - (E^1 I^1 \omega_{\z\z})(t,l-)
   &= - (E^2 I^2 \omega_{\z\z})(t,l+),
   &&t \geq 0
   \\
  (E^1 I^1 \omega_{\z\z})_\z(t,l-)
   &= (E^2 I^2 \omega_{\z\z})_\z(t,l+),
   &&t \geq 0
   \\
  (E^2 I^2 \omega_{\z\z})(t,1)
   &= 0,
   &&t \geq 0
   \\
  (E^2 I^2 \omega_{\z\z})_\z (t,1)
   &= \kappa \omega_t(t,1),
   &&t \geq 0,
   \\
  \omega(0,\z)
   &= \omega_0(\z),
   &&\z \in (0,1) \setminus \{l\}
   \\
  \omega_t(0,\z)
   &= \omega_1(\z),
   &&\z \in (0,1) \setminus \{l\}.
 \end{align*}
}For first reading, the reader may always have this special case in mind since it already includes most of the relevant features of a chain of Euler-Bernoulli beams.
The demonstration of the results in \cite{ChenEtAl_1987} is based on an energy multiplier method which provides a Lyapunov function for the Euler-Bernoulli beam system.
For example, in this case uniform exponential decay of the energy of the coupled system
 \[
  H(t)
   := \frac{1}{2} \sum_{j = 1}^m \int_{l^{j-1}}^{l^j} \rho^j \abs{\omega(t,\z)}^2 + E^j I^j \abs{\omega_{\z\z}(t,\z)}^2 \dd \z
   \leq M \ee^{\eta t} H(0),
   \quad
   t \geq 0
 \]
for some constants $M \geq 1$ and $\eta < 0$ which are independent of the initial data,
has been shown in \cite{ChenEtAl_1987} for a strictly positive damping parameter $\kappa > 0$ under the following additional structural constraints:
 \begin{enumerate}
  \item
   On each interval $(l^{j-1}, l^j)$, the mass density times cross sectional area $\rho(\z) = \rho^j$, the modulus of elasticity $E(\z) = E^j$ and the area moment of the cross section $I(\z) = I^j$ are constant.
  \item
   The parameters $\rho^j > 0$, $E^j > 0$ and $I^j > 0$ satisfy the monotonicity constraints
    \[
     \rho^j \leq \rho^{j+1},
      \quad
     E^j I^j \geq E^{j+1} I^{j+1},
      \quad
      j = 1, \ldots, m.
    \] 
 \end{enumerate}
In this paper, we are going to remove the first of these constraints, i.e.\ we  show the same uniform exponential stability result for arbitrary piecewise Lipschitz-continuous and strictly positive $\rho \in \Lip((l^{j-1}, l^j);\R)$ and $EI \in \Lip((l^{j-1}, l^j);\R)$ (note that this implies that for each junction point $l^j$, the one-sided limits $\rho(l^j-)$ and $\rho(l^j+)$ etc.\ exist) replacing the constant parameters $\rho^j$ and $E^j I^j$, but still satisfying a jump conditions for possible discontinuities of $\rho$ and $EI$ at the junction points:
 \begin{assumption}[Regularity of physical parameters and jump conditions]
  For the stability results will assume the following:
   \begin{align}
    \rho, EI &\in \Lip(l^{j-1},l^j)
    \quad \text{and uniformly positive},
    &&j = 1, \ldots, m
    \tag{{\bf R}}
    \label{R},
    \\
    \rho(l^j-) &\leq \rho(l^j+),
    \quad
    (EI)(l^j-) \geq (EI)(l^j+),
    &&j = 1, \ldots, m.
    \tag{{\bf M}}
    \label{M}
   \end{align}
 \end{assumption}

We prove our results within the framework of $C_0$-semigroups, applying a special case (for compact resolvents) of the Arendt-Batty-Lyubich-V\~u Theorem for asymptotic stability and the Gearhart-Pr\"uss-Huang Theorem for uniform exponential stability of $C_0$-semigroups.
Moreover, we consider a slight generalisation of the situation described above by allowing dynamic boundary feedback via impedance passive finite-dimensional control systems and which all are internally stable.
For well-posedness (here, in the sense that dissipativity of the interconnected system implies well-posedness with non-increasing energy for the system), we use abstract well-posedness results for so-called \emph{infinite-dimensional linear port-Hamiltonian systems} \cite{LeGorrecZwartMaschke_2004}, \cite{LeGorrecZwartMaschke_2005}, \cite{Villegas_2007}, \cite{AugnerJacob_2014} (which rely on the Lumer-Phillips Theorem), and also employ the techniques used in \cite{AugnerJacob_2014} for uniform exponential stabilisation of a (single) Euler-Bernoulli beam within the port-Hamiltonian framework to show the uniform exponential energy decay.
Note that in \cite{Augner_2020} quite general interconnection structures of infinite-dimensional port-Hamiltonian type have been considered, especially well-posedness and stability properties. The general setup considered there, however, is not enough to cover exponential stability of serially connected Euler-Bernoulli beams except for some very restrictive conditions on the dissipative structure of the interconnection.

This paper is organised as follows:
In Section \ref{well-posedness} we formally consider possible interconnection and boundary conditions leading to a dissipative system of joint Euler-Bernoulli beams.
More precisely, we give classes of boundary control and observation maps leading to an open loop impedance passive system, thus leading to a dissipative system for dissipative (linear) closure relations.
Using the abstract theory of infinite-dimensional linear port-Hamiltonian systems, this immediately implies well-posedness results in the sense that for any sufficiently regular initial data there is a unique solution with non-increasing energy and the solution depends continuously on the initial data.
In other words: The operator $\mathcal{A}$ governing the dynamics of the beam-observer-feedback-actuator system generates a strongly continuous contraction semigroup on a suitable energy state space $\mathcal{X}$.
Then, Section \ref{stability} is devoted to the discussion of stability properties.
We give sufficient conditions on the interconnection structure by means of dissipative static feedback or feedback via interconnection with an internally stable impedance passive finite-dimensional linear controller.
Here, the main results of that section and this manuscript are Theorem \ref{thm:asymptotic_stability} on asymptotic, i.e.\ strong, stability and Theorem \ref{thm:exp_stability} on uniform exponential stability, i.e.\ uniform exponential energy decay.
In Theorem \ref{exa:EB} we reformulate the previous well-posedness and stability results in the language of Euler-Bernoulli beam equations and show that our results cover the inhomogeneous beam versions of the uniform exponential stability results already presented in \cite{ChenEtAl_1987}, especially including a discussion of several relevant conservative boundary conditions on the non-dissipative end of the serial chain of Euler-Bernoulli beams, which already had been mentioned in \cite{ChenEtAl_1987}, but under the standing assumption of piecewise constant parameters $\rho$ and $EI$.
We conclude the paper with some final remarks in Section \ref{conclusion}.

\section{Well-posedness}
\label{well-posedness}

We start by discussing the well-posedness for systems of PDE modelling serially connected Euler-Bernoulli beams.
Slightly generalising the setup in \cite{ChenEtAl_1987}, we consider a decomposition $0 = l^0 < l^1 < \ldots < l^m = L$ of an interval $(0,L)$ ($L > 0$) and on each subinterval $(l^{j-1}, l^j)$ the Euler-Bernoulli beam equation
 \[
  \rho(\z) \omega_{tt}(t,\z) + (EI(\z) \omega_{\z\z}(t,\z))_{\z\z} =  0,
   \quad
   \z \in (l^{j-1}, l^j), \, j = 1, \ldots, m.
   \tag{EB}
   \label{EB}
 \]
where, in contrast to the situation in \cite{ChenEtAl_1987}, we allow for spatial dependence of $\rho$ and $EI$ on $\z \in (l^{j-1}, l^j)$. For the moment, it is enough to let $\rho, EI \in L_\infty(l^{j-1}, l^j)$ be uniformly positive, i.e.\ there is $\varepsilon > 0$ such that
 \[
  \rho(\z), EI(\z) \geq \varepsilon,
   \quad
   \text{a.e.\ } \z \in (l^{j-1}, l^j), \, j = 1, \ldots, m.
 \]
The total energy of this linear system is defined as the sum of kinetic energy and strain energy of all beams
 \begin{align*}
  H(t)
   &:= \frac{1}{2} \int_0^L \big( \rho(\z) \abs{\omega_t(t,\z)}^2 + EI(\z) \abs{\omega_{\z\z}(t,\z)}^2 \big) \dd \z
   \\
   &= \sum_{j=1}^m \frac{1}{2} \int_{l^{j-1}}^{l^j} \big( \rho(\z) \abs{\omega_t(t,\z)}^2 + EI(\z) \abs{\omega_{\z\z}(t,\z)}^2 \big) \dd \z
   =: \sum_{j=1}^m H^j(t).
 \end{align*}
For sufficiently regular solutions of the Euler-Bernoulli beam equations on each subinterval $(l^{j-1}, l^j)$, we formally obtain the power balance for the corresponding beam as
 \begin{align*}
  \frac{\dd}{\dd t} H^j(t)
   &= \Re \int_{l^{j-1}}^{l^j} \rho(\z) \omega_{tt}(t,\z) \overline{\omega_t(t,\z)} + EI(\z) \omega_{\z\z}(t,\z) \overline{\omega_{\z\z t}(t,\z)} \dd \z
   \\
   &= \Re \int_{l^{j-1}}^{l^j} - (EI \omega_{\z\z})_{\z\z}(t,\z) \overline{\omega_t(t,\z)} + EI(\z) \omega_{\z\z}(t,\z) \overline{\omega_{\z\z t}(t,\z)} \dd \z
   \\
   &= \Re \left[ - (EI \omega_{\z\z})_\z(t,l^j-) \overline{\omega_t(t,l^j-)} + (EI \omega_{\z\z})(t,l^j-) \overline{\omega_{t\z}(t,l^j-)} \right]
    \\
    &\quad
    - \Re \left[ - (EI \omega_{\z\z})_\z(t,l^{j-1}+) \overline{\omega_t(t,l^{j-1}+)} + (EI \omega_{\z\z})(t,l^{j-1}+) \overline{\omega_{t\z}(t,l^{j-1}+)} \right].
 \end{align*}
Putting these equations together, we obtain the change of total energy of sufficiently regular solutions as 
 \begin{align}
  \frac{\dd}{\dd t} H(t)
   = \sum_{j=1}^m \frac{\dd}{\dd t} H^j(t)
   &= \Re \left[ - (EI \omega_{\z\z})_\z(t,L-) \overline{\omega_t(t,L-)} + (EI \omega_{\z\z})(t,L-) \overline{\omega_{t\z}(t,L-)} \right]
    \nonumber \\
    &\quad
    - \Re \left[ - (EI \omega_{\z\z})_\z(t,0+) \overline{\omega_t(t,0+)} + (EI \omega_{\z\z})(t,0+) \overline{\omega_{t\z}(t,0+)} \right]
    \nonumber \\
    &\quad
    + \sum_{j=1}^{m-1} \Re \left[ - (EI \omega_{\z\z})_\z(t,l^j-) \overline{\omega_t(t,l^j-)} + (EI \omega_{\z\z})(t,l^j-) \overline{\omega_{t\z}(t,l^j-)} \right]
    \nonumber \\
    &\qquad
    - \Re \left[ - (EI \omega_{\z\z})_\z(t,l^j+) \overline{\omega_t(t,l^j+)} + (EI \omega_{\z\z})(t,l^j+) \overline{\omega_{t\z}(t,l^j+)} \right].
    \label{eqn:energy-balance}
 \end{align}
We see that, when aiming for interconnection of the beams in a dissipative way, the most natural way to do so is by imposing dissipative boundary conditions (which in this context may include conservative boundary conditions as well) at the left ($\z = l^0 = 0$) and right end ($\z = l^m = L$) and a dissipative interconnection at the junction points $\z = l^j$ ($j = 1, \ldots, m-1$).
To achieve the latter, at every junction point $l^j$ ($j = 1, \ldots, m-1$) we demand continuity conditions of the type
 \[
  \omega_t(t,l^j-)  = \omega_t(t,l^j +)
   \quad \text{ or } \quad
   - (EI \omega_{\z\z})_\z(t,l^j-) = - (EI \omega_{\z\z})_\z(t,l^j+)
 \]
and
 \[
  \omega_{t\z}(t,l^j-)  = \omega_{t\z}(t,l^j +)
   \quad \text{ or } \quad
   (EI \omega_{\z\z})(t,l^j-) = (EI \omega_{\z\z})(t,l^j+).
 \]
(In \cite{Pilkey_1969}, \cite{ChenEtAl_1987} it has been discussed that for dissipativity of the system, at least one of two state variables which are \emph{dual} (or \emph{complementary}) to each other has to be continuous. This justifies these conditions.)
We may thus distinguish between four different cases of static interconnections:
 \begin{enumerate}
  \item
   For $\omega_t(t,l^j-)  = \omega_t(t,l^j +)$ and $\omega_{t\z}(t,l^j-)  = \omega_{t\z}(t,l^j +)$:
   \[
    \left( \begin{array}{c} - (EI \omega_{\z\z})_\z(t,l^j-) + (EI \omega_{\z\z})_\z(t,l^j+) \\ (EI \omega_{\z\z})(t,l^j-) - (EI \omega_{\z\z})(t,l^j+) \end{array} \right) = - K^j \left( \begin{array}{c} \omega_t(t,l^j) \\ \omega_{t\z}(t,l^j) \end{array} \right)
   \]
  \item
   For $\omega_t(t,l^j-)  = \omega_t(t,l^j +)$ and $(EI \omega_{\z\z})(t,l^j-) = (EI \omega_{\z\z})(t,l^j+)$:
   \[
    \left( \begin{array}{c} - (EI \omega_{\z\z})_\z(t,l^j-) + (EI \omega_{\z\z})_\z(t,l^j+) \\ \omega_{t\z}(t,l^j-) - \omega_{t\z}(t, l^j+) \end{array} \right) = - K^j \left( \begin{array}{c} \omega_t(t,l^j) \\ (EI \omega_{\z\z})(t,l^j) \end{array} \right).
   \]
  \item
   For $- (EI \omega_{\z\z})_\z(t,l^j-) = - (EI \omega_{\z\z})_\z(t,l^j+)$ and $\omega_{t\z}(t,l^j-)  = \omega_{t\z}(t,l^j +)$:
   \[
    \left( \begin{array}{c} \omega_t(t,l^j-) - \omega_t(t,l^j+) \\ (EI \omega_{\z\z})(t,l^j-) - (EI \omega_{\z\z})(t,l^j+) \end{array} \right) = - K^j \left( \begin{array}{c} - (EI \omega_{\z\z})_\z(t,l^j) \\ \omega_{t\z}(t,l^j) \end{array} \right).
   \]
  \item
   For $- (EI \omega_{\z\z})_\z(t,l^j-) = - (EI \omega_{\z\z})_\z(t,l^j+)$ and $(EI \omega_{\z\z})(t,l^j-) = (EI \omega_{\z\z})(t,l^j+)$:
   \[
    \left( \begin{array}{c} \omega_t(t,l^j-) - \omega_t(t,l^j+) \\ \omega_{t\z}(t,l^j-) - \omega_{t\z}(t, l^j+) \end{array} \right) = - K^j \left( \begin{array}{c} - (EI \omega_{\z\z})_\z(t,l^j) \\ (EI \omega_{\z\z})(t,l^j) \end{array} \right).
   \]
 \end{enumerate}
In each case $K^j \in \K^{2 \times 2}$ denotes a matrix with positive semidefinite symmetric part $\Sym K = \frac{K + K^*}{2}$.
(At first reading the reader might consider the special case $m = 2$ and the particular interconnection condition
 \begin{align*}
  \omega_t(l^1+)
   &= \omega_t(t,l^1-)
   \\
  \omega_{t\z}(t,l^1+)
   &= \omega_{t\z}(t,l^1-)
   \\
  \left( \begin{array}{c} - (EI \omega_{\z\z})_\z(t,l^1-) + (EI \omega_{\z\z})_\z(t,l^1+) \\ (EI \omega_{\z\z})(t,l^1-) - (EI \omega_{\z\z})(t,l^1+) \end{array} \right)
   &= - K^1 \left( \begin{array}{c} \omega_t(t,l^1) \\ \omega_{t\z}(t,l^1) \end{array} \right)
 \end{align*}
to make the following statements easier digestible.)
Reformulating this problem in a more abstract way, makes the theory of infinite-dimensional linear port-Hamiltonian systems, cf.\ e.g.\ \cite{LeGorrecZwartMaschke_2005}, \cite{AugnerJacob_2014}, applicable, and we, therefore, may easily deduce well-posedness from abstract theory, simply by checking dissipativity.
First, we exploit that $\rho$ and $EI$ are bounded and uniformly positive and write
 \begin{align*}
  \tilde \H(\z)
   &:= \left[ \begin{array}{cc} \tilde \H_1(\z) & \\ & \tilde \H_2(\z) \end{array} \right]
   := \left[ \begin{array}{cc} \rho(\z)^{-1} & \\ & EI(\z) \end{array} \right]
   \in \K^{2 \times 2},
   \\
  \tilde x(t,\z)
   &:= \left( \begin{array}{c} \tilde x_1(t,\z) \\ \tilde x_2(t,\z) \end{array} \right)
   := \left( \begin{array}{c} \rho(\z) \omega_t(t,\z) \\ \omega_{\z\z}(t,\z) \end{array} \right)
   \in \K^2,
 \end{align*}
so that the Euler-Bernoulli beam equations may be equivalently expressed as
 \[
  \frac{\partial}{\partial t} \tilde x(t,\z)
   = \left[ \begin{array}{cc} 0 & -1 \\ 1 & 0 \end{array} \right] \frac{\partial^2}{\partial \z^2} (\tilde \H(\z) \tilde x(t,\z)),
   \quad
   t \geq 0, \, \z \in (l^{j-1}, l^j), \, j = 1, \ldots, m. 
 \]
With $\tilde P_2 = \left[ \begin{array}{cc} 0 & -1 \\ 1 & 0 \end{array} \right]$, this almost looks like an infinite-dimensional linear port-Hamiltonian system of order 2 as considered in \cite{LeGorrecZwartMaschke_2005} or \cite{AugnerJacob_2014}, but due to the possible discontinuities in the junction points it is not quite yet.
However, by performing a parameter transformation and writing the resulting PDE as a system of PDE on the unit interval $(0,1)$,
 \begin{align*}
  \H^j(\z) &:= \tilde \H((1-\z) l^{j-1} + \z l^j),
   \quad
   x^j(t,\z) := \tilde x(t,(1-\z) l^{j-1} + \z l^j),
   \quad
   \z \in (0,1), \, j = 1, \ldots, m
   \\
  \H(\z) &:= \left[ \begin{array}{ccc} \H^1(\z) && \\ & \ddots & \\ && \H^m(\z) \end{array} \right]
   \in \K^{2m \times 2m}
   \\
  x(t,\z) &:= \left( \begin{array}{c} x^1(t,\z) \\ \vdots \\ x^m(t,\z) \end{array} \right)
   \in \K^{2m}
   \\
   P_2 &= \left[ \begin{array}{ccc} \frac{1}{(l^1 - l^0)^2} \left[ \begin{array}{cc} 0 & -1 \\ 1 & 0 \end{array} \right] && \\ & \ddots & \\ && \frac{1}{(l^m - l^{m-1})^2} \left[ \begin{array}{cc} 0 & -1 \\ 1 & 0 \end{array} \right] \end{array} \right]
   \in \K^{2m \times 2m},
   \\
  P_1
   &= P_0 = 0 \in \K^{2m \times 2m}
 \end{align*}
we find that the Euler-Bernoulli equations take the port-Hamiltonian form
 \[
  \frac{\partial}{\partial t} x(t,\z)
   = \left( P_2 \frac{\partial^2}{\partial \z^2} \right) (\H(\z) x(t,\z))
   =: \left( \mathfrak{A} x(t) \right)(\z),
   \quad
   t \geq 0, \, \z \in (0,1).
 \]
Let the \emph{energy state space} $X= L_2(0,1;\K^{2m})$ be equipped with the energy inner product
 \[
  \sp{x}{y}_X
   := \int_0^1 \sp{x(\z)}{\H(\z) y(\z)}_{\K^{2m}} \dd \z,
   \quad
   x, y \in X.
 \]
By \cite[Theorem 4.1]{LeGorrecZwartMaschke_2005}, the operator $A = \mathfrak{A}|_{\dom(A)}$ defined as the restriction of $\mathfrak{A}$ to any domain $\dom(A) \subseteq \dom(\mathfrak{A}) = \{ x \in L_2(0,1;\K^{2m}): \, \H x \in H^2(0,1;\K^{2m}) \}$ of the form \[ \dom(A) = \{ x \in \dom(\mathfrak{A}): \, W \left( \begin{array}{c} (\H x)(0) \\ (\H x)(1) \end{array} \right) = 0 \} \quad \text{for some full rank matrix } W \in \K^{2m \times 4m} \] generates a contractive $C_0$-semigroup on $X$ if and only if $A$ is dissipative. I.e., the boundary (and here: interconnection) conditions restricting $\dom(\mathfrak{A})$ to a linear subspace $\dom(A)$ need to ensure that
 \[
  \Re \sp{Ax}{x}_X
   \leq 0,
   \quad
   x \in \dom(A).
 \]
Note that thanks to the uniform boundedness and positivity of $\rho$ and $EI$, $\sp{\cdot}{\cdot}_X$ induces a norm on $X$ which is equivalent to the standard inner product on $L_2(0,1;\K^{2m})$.
A similar statement holds true whenever the system of beams is interconnected by a finite dimensional linear control system $\Sigma_c = (A_c, B_c, C_c, D_c)$, i.e.\
 \begin{align*}
  \frac{\dd}{\dd t} x_c(t)
   &= A_c x_c(t) + B_c u_c(t),
   \quad
   t \geq 0
   \\
  y_c(t)
   &= C_c x_c(t) + D_c u_c(t),
   \quad
   t \geq 0.
 \end{align*}
Here, $x_c(t)$ lies in the controller state space $X_c$, a finite dimensional Hilbert space, and the Hilbert spaces $U_c$ and $Y_c$ are the finite dimensional control and observation spaces for the control system.
To formulate this well-posedness result rigorously, we introduce the following boundary control and observation operators for the Euler-Bernoulli beam system in port-Hamiltonian form.

\begin{dfntn}[Pointwise Control and Observation Operators]
Let matrices $W_B^j, W_C^j \in \K^{2 \times 4}$ be given such that $\left[ \begin{smallmatrix} W_B^j \\ W_C^j \end{smallmatrix} \right] \in \K^{4 \times 4}$ is invertible.
We define the linear operators $\mathfrak{B}^0, \mathfrak{C}^0, \mathfrak{B}^m, \mathfrak{C}^m: \dom(\mathfrak{A}) \rightarrow \K^2$ by 
 \begin{align*}
  \left( \begin{array}{c} \mathfrak{B}^0 x \\ \mathfrak{C}^0 x \end{array} \right)
   &:= \left[ \begin{array}{c} W_B^0 \\ W_C^0 \end{array} \right] \left( \begin{array}{c} (\H^1 x^1)(0) \\ \tilde P_2 (\H^1 x^1)'(0) \end{array} \right)
   = \left[ \begin{array}{c} W_B^0 \\ W_C^0 \end{array} \right] \left( \begin{array}{c} \omega_t(t,0) \\ (EI \omega_{\z\z})(t,0) \\ - (EI \omega_{\z\z})_\z(t,0) \\ \omega_{t\z}(t,0) \end{array} \right)
   \\
  \left( \begin{array}{c} \mathfrak{B}^m x \\ \mathfrak{C}^m x \end{array} \right)
   &:= \left[ \begin{array}{c} W_B^m \\ W_C^m \end{array} \right] \left( \begin{array}{c} (\H^m x^m)(1) \\ - \tilde P_2 (\H^m x^m)'(1) \end{array} \right)
   = \left[ \begin{array}{c} W_B^m \\ W_C^m \end{array} \right] \left( \begin{array}{c} \omega_t(t,L) \\ (EI \omega_{\z\z})(t,L) \\ (EI \omega_{\z\z})_\z(t,L) \\ - \omega_{t\z}(t,L) \end{array} \right)
 \end{align*}
and for each junction $j \in \{1, \ldots, m-1\}$, depending the chosen type of continuity conditions, we define linear maps $\mathfrak{B}_0^j$, $\mathfrak{B}^j$ and $\mathfrak{C}^j: \dom(\mathfrak{A}) \rightarrow \K^2$ as follows:
 \begin{enumerate}
  \item
   For the case where $\omega_t(t,l^j-)  = \omega_t(t,l^j +)$ and $\omega_{t\z}(t,l^j-)  = \omega_{t\z}(t,l^j +)$ are continuous in $l^j$:
 \begin{align*}
  \mathfrak{B}_0^j x
   &:=  \left( \begin{array}{c} (\H^j_1 x^j_1)(1) - (\H^{j+1}_1 x^{j+1}_1)(0) \\ (\H^j_1 x^j_1)'(1) - (\H^{j+1}_1 x^{j+1}_1)'(0) \end{array} \right)
   = \left( \begin{array}{c} \omega_t(t,l^j-) - \omega_t(t,l^j+) \\ \omega_{t\z}(t,l^j-) - \omega_{t\z}(t,l^j+) \end{array} \right)
   \\
  \left( \begin{array}{c} \mathfrak{B}^j x \\ \mathfrak{C}^j x \end{array} \right)
   &:= \left( \begin{array}{c} - (\H^j_2 x^j_2)'(1) + (\H^{j+1}_2 x^{j+1}_2)'(0) \\ (\H^j_2 x^j_2)(1) - (\H^{j+1}_2 x^{j+1}_2)(0) \\ \tfrac{1}{2} ((\H^j_1 x^j_1)(1) + (\H^{j+1}_1 x^{j+1}_1)(0)) \\ \tfrac{1}{2} ((\H^j_1 x^j_1)'(1) + (\H^{j+1}_1 x^{j+1}_1)'(0)) \end{array} \right)
   = \left( \begin{array}{c} - (EI \omega_{\z\z})_\z(t,l^j-) + (EI \omega_{\z\z})_\z(t,l^j+) \\ (EI \omega_{\z\z})(t,l^j-) - (EI \omega_{\z\z})(t,l^j+) \\ \omega_t(t,l^j) \\ \omega_{t\z}(t,l^j) \end{array} \right)
 \end{align*}
  where for any spatial dependent quantity $f$ we write $f(\z) := f(\z+) = f(\z-)$ whenever the two one-sided limits exist and coincide.
  \item
   For the case $\omega_t(t,l^j-)  = \omega_t(t,l^j+)$ and $(EI \omega_{\z\z})(t,l^j-) = (EI \omega_{\z\z})(t,l^j+)$:
 \begin{align*}
  \mathfrak{B}_0^j x
   &:=  \left( \begin{array}{c} (\H^j_1 x^j_1)(1) - (\H^{j+1}_1 x^{j+1}_1)(0) \\ (\H^j_2 x^j_2)(1) - (\H^{j+1}_2 x^{j+1}_2)(0) \end{array} \right)
   = \left( \begin{array}{c} \omega_t(t,l^j-) - \omega_t(t,l^j+) \\ (EI \omega_{\z\z})(t,l^j-) - (EI \omega_{\z\z})(t,\z^j+) \end{array} \right)
   \\
  \left( \begin{array}{c} \mathfrak{B}^j x \\ \mathfrak{C}^j x \end{array} \right)
   &:= \left( \begin{array}{c} - (\H^j_2 x^j_2)'(1) + (\H^{j+1}_2 x^{j+1}_2)'(0) \\ (\H^j_1 x^j_1)'(1) - (\H^{j+1}_1 x^{j+1}_1)'(0) \\ \tfrac{1}{2} ((\H^j_1 x^j_1)(1) + (\H^{j+1}_1 x^{j+1}_1)(0)) \\ \tfrac{1}{2} ((\H^j_2 x^j_2)(1) + (\H^{j+1}_2 x^{j+1}_2)(0)) \end{array} \right)
   = \left( \begin{array}{c} - (EI \omega_{\z\z})_\z(t,l^j-) + (EI \omega_{\z\z})_\z(t,l^j+) \\ \omega_{t\z}(t,l^j-) - \omega_{t\z}(t,l^j+) \\ \omega_t(t,l^j) \\ (EI \omega_{\z\z})(t,l^j)  \end{array} \right)
 \end{align*}
  \item
   For $- (EI \omega_{\z\z})_\z(t,l^j-) = - (EI \omega_{\z\z})_\z(t,l^j+)$ and $\omega_{t\z}(t,l^j-)  = \omega_{t\z}(t,l^j +)$:
 \begin{align*}
  \mathfrak{B}_0^j x
   &:=  \left( \begin{array}{c} - (\H^j_2 x^j_2)'(1) - (\H^{j+1}_2 x^{j+1}_2)'(0) \\ (\H^j_1 x^j_1)'(1) - (\H^{j+1}_1 x^{j+1}_1)'(0) \end{array} \right)
   = \left( \begin{array}{c} - (EI \omega_{\z\z})_\z(t,l^j-) - (EI \omega_{\z\z})_\z(t,l^j+) \\ \omega_{t\z}(t,l^j-) - \omega_{t\z}(t,l^j+) \end{array} \right)
   \\
  \left( \begin{array}{c} \mathfrak{B}^j x \\ \mathfrak{C}^j x \end{array} \right)
   &:= \left( \begin{array}{c} - (\H^j_1 x^j_1)(1) + (\H^{j+1}_1 x^{j+1}_1)(0) \\ (\H^j_2 x^j_2)(1) - (\H^{j+1}_2 x^{j+1}_2)(0) \\ - \tfrac{1}{2} ((\H^j_2 x^j_2)'(1) + (\H^{j+1}_2 x^{j+1}_2)'(0)) \\ \tfrac{1}{2} ((\H^j_1 x^j_1)'(1) + (\H^{j+1}_1 x^{j+1}_1)'(0)) \end{array} \right)
   = \left( \begin{array}{c} \omega_t(t,l^j-) - \omega_t(t, l^j +) \\ (EI \omega_{\z\z})(t,l^j-) - (EI \omega_{\z\z})(t,l^j+) \\ - (EI \omega_{\z\z})_\z(t, l^j) \\ \omega_{t\z}(t,l^j) \end{array} \right)
 \end{align*}
  \item
   For $- (EI \omega_{\z\z})_\z(t,l^j-) = - (EI \omega_{\z\z})_\z(t,l^j+)$ and $(EI \omega_{\z\z})(t,l^j-) = (EI \omega_{\z\z})(t,l^j+)$:
 \begin{align*}
  \mathfrak{B}_0^j x
   &:=  \left( \begin{array}{c} - (\H^j_2 x^j_2)'(1) - (\H^{j+1}_2 x^{j+1}_2)'(0) \\ (\H^j_2 x^j_2)(1) - (\H^{j+1}_2 x^{j+1}_2)(0) \end{array} \right)
   = \left( \begin{array}{c} - (EI \omega_{\z\z})_\z(t, l^j-) - (EI \omega_{\z\z})_\z(t, l^j+) \\ (EI \omega_{\z\z})(t, l^j-) - (EI \omega_{\z\z})(t, l^j+) \end{array} \right)
   \\
  \left( \begin{array}{c} \mathfrak{B}^j x \\ \mathfrak{C}^j x \end{array} \right)
   &:= \left( \begin{array}{c} (\H^j_1 x^j_1)(1) - (\H^{j+1}_1 x^{j+1}_1)(0) \\ (\H^j_1 x^j_1)'(1) - (\H^{j+1}_1 x^{j+1}_1)'(0) \\ - \tfrac{1}{2} ((\H^j_2 x^j_2)'(1) + (\H^{j+1}_2 x^{j+1}_2)'(0)) \\ \tfrac{1}{2} ((\H^j_2 x^j_2)(1) + (\H^{j+1}_2 x^{j+1}_2)(0)) \end{array} \right)
   = \left( \begin{array}{c} \omega_t(t,l^j-) - \omega_t(t,l^j+) \\ \omega_{t\z}(t,l^j-) - \omega_{t\z}(t, l^j+) \\ - (EI \omega_{\z\z})_\z(t,l^j) \\ (EI \omega_{\z\z})(t,l^j) \end{array} \right).
 \end{align*}
 \end{enumerate}
\end{dfntn}

When interconnecting with a control system, or closing the system via dissipative boundary feedback and interconnection at the junction points, it is convenient to have an \emph{impedance passive} system.
We define the linear operators $\mathfrak{A}_0$, $\mathfrak{B}$ and $\mathfrak{C}$ by
 \begin{align*}
  \mathfrak{A}_0
   &= \mathfrak{A}|_{\ker \mathfrak{B}_0}
   = \mathfrak{A}|_{\cap_{j=1}^{m-1} \ker \mathfrak{B}_0^j}:
   \quad
   \dom(\mathfrak{A}_0) \subseteq X \rightarrow X
   \\
  \mathfrak{B}
   &= (\mathfrak{B}^j)_{j=0}^m:
   \quad
   \dom(\mathfrak{A}_0) \subseteq X \rightarrow U := \K^{2(m+1)}
   \\
  \mathfrak{C}
   &= (\mathfrak{C}^j)_{j=0}^m:
   \quad
   \dom(\mathfrak{A}_0) \subseteq X \rightarrow Y := \K^{2(m+1)}.
 \end{align*}
For this choice, the power balance \eqref{eqn:energy-balance} then implies that for sufficiently regular solutions of $x_t(t,\z) = \mathfrak{A}_0 x(t,\z)$ the energy changes as $\frac{\dd}{\dd t} \frac{1}{2} \norm{x}_X^2 = \Re \sp{\mathfrak{A}_0 x}{x}_X$, where
 \begin{align*}
  \Re \sp{\mathfrak{A}_0 x}{x}_X
   &= \sum_{j = 1}^{m-1} \Re \sp{\mathfrak{B}^j x}{\mathfrak{C}^j x}
    + \Re \sp{- \tilde P_2 (\H^1 x^1)'(0)}{(\H^1 x^1)(0)}
    \\
    &\quad
    + \Re \sp{\tilde P_2 (\H^m x^m)'(1)}{(\H^m x^m)(1)} ,
   \quad
   x \in \dom(\mathfrak{A}^0).
 \end{align*}
From here, conditions on the matrices $\left[ \begin{array}{c} W_B^0 \\ W_C^0 \end{array} \right]$ and $\left[ \begin{array}{c} W_B^m \\ W_C^m \end{array} \right]$ may lead to impedance passivity of the system $\mathfrak{S}_0 = (\mathfrak{A}_0, \mathfrak{B}, \mathfrak{C})$.

\begin{lmm}
The triplet $(\mathfrak{A}_0, \mathfrak{B}, \mathfrak{C})$ is impedance passive, i.e.\
 \[
  \Re \sp{\mathfrak{A}_0 x}{x}
   \leq \Re \sp{\mathfrak{B} x}{\mathfrak{C} x}_U
   = \sum_{j=0}^m \Re \sp{\mathfrak{B}^j x}{\mathfrak{C}^j x}_{\K^2},
   \quad
   x \in \dom(\mathfrak{A}_0),
 \]
if and only if
 \[
  \left[ \begin{array}{c} W_B^0 \\ W_C^0 \end{array} \right]^{\ast} \left[ \begin{array}{cc} 0 & I \\ I & 0 \end{array} \right] \left[ \begin{array}{c} W_B^0 \\ W_C^0 \end{array} \right] - \left[ \begin{array}{cc} 0 & I \\ I & 0 \end{array} \right] \geq 0,
   \quad
  \left[ \begin{array}{c} W_B^m \\ W_C^m \end{array} \right]^{\ast} \left[ \begin{array}{cc} 0 & 1 \\ 1 & 0 \end{array} \right] \left[ \begin{array}{c} W_B^m \\ W_C^m \end{array} \right] - \left[ \begin{array}{cc} 0 & I \\ I & 0 \end{array} \right] \geq 0
 \]
are both positive semidefinite.
\end{lmm}

\textbf{Proof.}
This can easily be proved using the energy balance \eqref{eqn:energy-balance}.
\qed

\begin{assumption}[Impedance-passivity]
The triplet $\mathfrak{S}_0 = (\mathfrak{A}_0, \mathfrak{B}, \mathfrak{C})$ is impedance passive.
\end{assumption}

The class of control systems we consider for linear closing of the open loop system $(\mathfrak{A}_0, \mathfrak{B}, \mathfrak{C})$ is always assumed to be a linear finite dimensional controller. 

\begin{assumption}[Finite-dimensional, linear control system]
Let $X_c$ be any finite dimensional $\K$-Hilbert space (including the possible choice $X_c = \{0\}$ leading to static boundary feedback) and (w.l.o.g.) let $U_c = Y_c = \K^{2(m+1)}$ and $A_c \in \B(X_c)$, $B_c \in \B(U_c,X_c)$, $C_c \in \B(X_c, Y_c)$ and $D_c \in \B(U_c,Y_c)$ be bounded linear operators, defining the finite dimensional linear control system $\Sigma_c = (A_c, B_c, C_c, D_c)$ with dynamics
 \begin{align*}
   \frac{\dd}{\dd t} x_c(t) 
    &= A_c x_c(t) + B_c u_c(t)
    \\
   y_c(t)
    &= C_c x_c(t) + D_c u_c(t),
    \quad
    t \geq 0.
 \end{align*}
\end{assumption}

We consider the standard feedback interconnection between the chain of Euler-Bernoulli beams and the controller which is given by
 \[
  \mathfrak{B} x(t)
   = - y_c(t),
   \quad
  u_c(t)
   = \mathfrak{C} x(t),
   \qquad
   t \geq 0.
 \]
The dynamics of the interconnected system is then described by the abstract Cauchy problem
 \begin{equation}
  \begin{cases}
   \frac{\dd}{\dd t} (x,x_c)(t)
    = \mathcal{A} (x,x_c)(t),
    &
    t \geq 0
    \\
   (x,x_c)(0)
    = (x_0, x_{c,0})
   \end{cases}
  \tag{ACP}
  \label{eqn:ACP}
 \end{equation}
for some given initial data $(x_0,x_{c,0}) \in \mathcal{X} := X \times X_c$. Here, the product Hilbert space $\mathcal{X}$ is equipped with the inner product
 \[
  \sp{(x,x_c)}{(z,z_c)}_{\mathcal{X}}
   = \sp{x}{z}_X + \sp{x_c}{z_c}_{X_c},
    \quad
    (x,x_c), (z,z_c) \in \mathcal{X} = X \times X_c
 \]
and the linear operator $\mathcal{A}: \dom(\mathcal{A}) \subseteq \mathcal{X} \rightarrow \mathcal{X}$ is defined by
 \begin{align*}
  \mathcal{A} \left( \begin{array}{c} x \\ x_c \end{array} \right)
   &= \left[ \begin{array}{cc} \mathfrak{A} & 0 \\ B_c \mathfrak{C} & A_c \end{array} \right] \left( \begin{array}{c} x \\ x_c \end{array} \right)
   \\
  \dom(\mathcal{A})
   &= \left\{ (x,x_c) \in \dom(\mathfrak{A}) \times X_c: \quad \mathfrak{B}_0 x = 0, \, \mathfrak{B} x = - (C_c x_c + D_c \mathfrak{C} x) \right\}
   \\
   &= \left\{ (x,x_c) \in \dom(\mathfrak{A}_0) \times X_c: \quad \mathfrak{B} x = - (C_c x_c + D_c \mathfrak{C} x) \right\}
 \end{align*}
where $\mathfrak{B}_0: \dom(\mathfrak{B}_0) = \dom(\mathfrak{A}) \subseteq X \rightarrow \K^{2(m-1)}$ is defined by
 \[
  \mathfrak{B}_0 x
   = (\mathfrak{B}_0^j x^j)_{j=1}^{m-1}.
 \]
 
 \begin{prpstn}[Well-posedness of the abstract Cauchy problem]
  The operator $\mathcal{A}$ generates a contractive $C_0$-semigroup $(\mathcal{T}(t))_{t \geq 0}$ on $\mathcal{X}$ if and only if $\mathcal{A}$ is dissipative, i.e.\
   \[
    \Re \sp{\mathcal{A} (x,x_c)}{(x,x_c)}_{\mathcal{X}}
     \leq 0,
     \quad
     (x,x_c) \in \dom(\mathcal{A}) \times X_c.
   \]
  Moreover, in that case the operator $\mathcal{A}$ has compact resolvent.
  In particular, this is the case if $\mathfrak{S}_0 = (\mathfrak{A}_0, \mathfrak{B}, \mathfrak{C})$ and $\Sigma_c = (A_c, B_c, C_c, D_c)$ are impedance passive, i.e.\
   \begin{align*}
    \Re \sp{\mathfrak{A}_0 x}{x}_X
     &\leq \Re \sp{\mathfrak{B} x}{\mathfrak{C} x}_U,
     &&x \in \dom(\mathfrak{A}_0)
     \\
    \Re \sp{A_c x_c + B_c u_c}{x_c}_{X_c}
     &\leq \Re \sp{C_c x_c + D_c u_c}{u_c}_{U_c},
     &&x_c \in x_c, \, u_c \in U_c.
   \end{align*}
 \end{prpstn}
 
\textbf{Proof.}
See \cite[Theorem 3.1]{AugnerJacob_2014}.\qed

In other words, for every linear closure via static feedback or a dynamic linear, finite dimensional control system such that the resulting interconnected system is dissipative, and for any initial data $(x_0, x_{c,0}) \in X \times X_c$ there is a unique strong solution $(x,x_c) \in C(\R_+; X \times X_c)$ of the abstract Cauchy problem with non-increasing energy
 \[
  H^{\mathrm{tot}}(t)
   = \frac{1}{2} \norm{x(t)}_X^2 + \frac{1}{2} \norm{x_c(t)}_{X_c}^2,
   \quad
   t \geq 0.
 \]

\section{Stability properties}
\label{stability}

\begin{figure}
	\centering
	\includegraphics{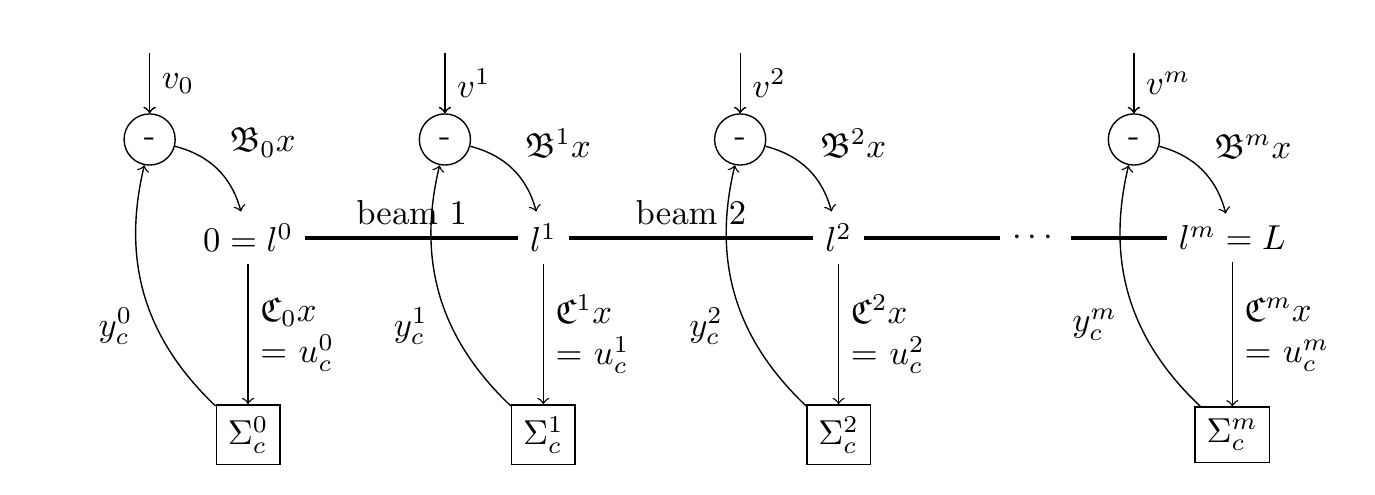}
	\caption{At each junction point $l^j$ the beam ends meeting there are interconnected via a finite dimensional control system $\Sigma_c^j$.}
	\label{img:interconnection_structure}
\end{figure}

Having established well-posedness for every set of dissipative interconnection and boundary conditions, in this section we investigate stability properties under a slightly more restrictive condition on the finite dimensional control system.
Namely, we assume that $\Sigma_c$ has block diagonal form in the following sense: There are control systems $\Sigma_c^j = (A_c^j, B_c^j, C_c^j, D_c^j)$, $j = 0, 1, \ldots, m$, defined on finite-dimensional controller state spaces $X_c^j$ such that $X_c = \prod_{i=0}^m X_c^j$ and input and output spaces $U_c^j = Y_c^j = \K^2$ such that $U_c = \prod_{i=0}^m U_c^j = \prod_{i=0}^m Y_c^j = Y_c$, and linear operators $A_c^j \in \B(X_c^j)$, $B_c^j \in \B(U_c^j, X_c^j)$, $C_c^j \in \B(X_c^j; Y_c^j)$ and $D_c^j \in \B(U_c^j; Y_c^j)$ such that $A_c = \operatorname{diag}(A_c^j)_{j=0,\ldots,m}$, $B_c = \operatorname{diag}(B_c^j)_{j=0,\ldots,m}$ etc.\ are block-diagonal linear operators.
The standard feedback interconnection then reads as
 \[
  \mathfrak{B}^j x = - y_c^j,
   \quad
  \mathfrak{C}^j x = u_c^j,
   \quad
   j = 0, 1, \ldots, m
 \]
with the dynamics of the finite dimensional controllers being modelled by
 \begin{align*}
   \frac{\dd}{\dd t} x_c^j(t) 
    &= A_c^j x_c^j(t) + B_c u_c^j(t)
    \\
   y_c^j(t)
    &= C_c^j x_c^j(t) + D_c^j u_c^j(t),
    \quad
    j = 0, 1, \ldots, m, \,
    t \geq 0,
 \end{align*}
see Figure \ref{img:interconnection_structure}.
 
\begin{assumption}[Passive and internally stable controllers]
\label{assmpt:passive_stable_controller}
 All control systems $\Sigma_c^j = (A_c^j, B_c^j, C_c^j, D_c^j)$ are impedance passive and for $j = 0, 1, \ldots, m$ the following conditions hold:
  \begin{enumerate}
   \item
    \label{assmpt:passive_stable_controller-i}
    There is $\kappa_j > 0$ such that
     \[
      \Re \sp{A_c^j x_c^j + B_c^j u_c^j}{x_c^j}_{X_c^j}
       \leq \Re \sp{C_c^j x_c^j + D_c^j u_c^j}{u_c^j}_{U_c^j} - \kappa_j \abs{D_c^j u_c^j},
       \quad
       x_c^j \in X_c^j, \, u_c^j \in U_c^j.
     \]
    \item
     $\ker D_c^j \subseteq \ker B_c^j$.
    \item
     $A_c^j$ has spectrum $\sigma(A_c^j) \subseteq \C_0^- := \{\lambda \in \C: \Re \lambda < 0\}$, i.e.\ the semigroup $(\ee^{t A_c^j})_{t \geq 0}$ is exponentially stable.
  \end{enumerate}
\end{assumption}

\begin{rmrk}[Collocated input and output]
 If $B_c^j = (C_c^j)^\ast$ for all $j = 0, 1, \ldots, m$, i.e.\ collocated input and output, and each $D_c^j$ is either a diagonal matrix with non-negative entries or the symmetric part $\Sym(D_c^j)$ is positive definite, then Assumption \ref{assmpt:passive_stable_controller}.\eqref{assmpt:passive_stable_controller-i} is satisfied.
 Namely, impedance passivity implies that each $A_c^j$ is dissipative, and in either case one has $\Re \sp{D_c^j u_c^j}{u_c^j}_{U_c^j} = \abs{(\Sym D_c^j)^{1/2} u_c^j}_{U_c^j}^2 \geq \kappa \abs{D_c^j u_c^j}^2$ for some $\kappa > 0$.
\end{rmrk}

\begin{rmrk}
\label{rem:on_assmpt_3.1}
 Using standard feedback interconnection between the impedance passive systems $\Sigma_c$ and $\mathfrak{S}$ in the definition of the operator $\mathcal{A}$, in connection with Assumption \ref{assmpt:passive_stable_controller} gives that
  \begin{align*}
   \Re \sp{\mathcal{A} (x,x_c)}{(x,x_c)}_{\mathcal{X}}
    &= \Re \sp{\mathfrak{A} x}{x}_X + \Re \sp{A_c x_c + B_c \mathfrak{C} x}{- \mathfrak{B} x}_{X_c}
    \\
    &\leq \Re \sp{\mathfrak{B} x}{\mathfrak{C} x}_{U} + \Re \sp{\mathfrak{C} x}{- \mathfrak{B} x}_U - \sum_{j=0}^m \kappa_j \abs{D_c^j \mathfrak{C}^j x}^2
    \\
    &= - \sum_{j=0}^m \kappa_j \abs{D_c^j \mathfrak{C}^j x}^2,
    \quad
    (x,x_c) \in \dom(\mathcal{A}).
  \end{align*}
\end{rmrk}

Under this structural assumption and some slight regularity conditions on $\rho$ and $EI$, we are able to formulate the following asymptotic stability result for linear boundary damping feedback at one of the two ends of the chain of Euler-Bernoulli beams.

\begin{thrm}[Asymptotic Stability]
 \label{thm:asymptotic_stability}
 Assume that $\rho$ and $EI$ have regularity \eqref{R}, the control systems satisfy Assumption \ref{assmpt:passive_stable_controller} and for the operator $\mathcal{A}$ defined in Section \ref{well-posedness} it holds
  \begin{equation}
   \Re \sp{\mathcal{A} (x,x_c)}{(x,x_c)}_{\mathcal{X}}
    \leq - \kappa \abs{\mathfrak{R} x}_{\K^4}^2,
    \quad
    x \in \dom(\mathcal{A})
    \label{eqn:diss_asymptoctic}
  \end{equation}
 where $\kappa > 0$ and $\mathfrak{R}: \dom(\mathfrak{A}) \rightarrow \K^4$ is one of the functions
  \[
   \left( \begin{array}{c} (\H^1_1 x^1_1)(0) \\ (\H^1_1 x^1_1)'(0) \\ (\H^1_2 x^1_2)(0) \\ (\H^m_2 x^m_2)'(1) \end{array} \right) \, \text{or} \,
    \left( \begin{array}{c} (\H^1_1 x^1_1)(0) \\ (\H^1_1 x^1_1)'(0) \\ (\H^1_2 x^1_2)'(0) \\ (\H^m_2 x^m_2)(1) \end{array} \right) \, \text {or} \,
    \left( \begin{array}{c} (\H^1_1 x^1_1)(0) \\ (\H_2^1 x_2^1)(0) \\ (\H^1_2 x^1_2)'(0) \\ (\H^m_1 x^m_1)'(1) \end{array} \right) \, \text{or} \, 
    \left( \begin{array}{c} (\H^1_1 x^1_1)'(0) \\ (\H^1_2 x^1_2)(0) \\ (\H^1_2 x^1_2)'(0) \\ (\H^m_1 x^m_1)(1) \end{array} \right).
  \]
 Then the $C_0$-semigroup $(\mathcal{T}(t))_{t \geq 0}$ generated by $\mathcal{A}$ is asymptotically stable on $\mathcal{X}$, i.e.\ for every initial value $(x_0, x_{c,0}) \in \mathcal{X}$ one has
  \[
   \mathcal{T}(t) (x_0, x_{c,0}) \rightarrow 0.
  \]
 In particular, $\sigma(\mathcal{A}) = \sigma_p(\mathcal{A}) \subseteq \C_0^-$.
\end{thrm}

\textbf{Proof.}
Since the operator $\mathcal{A}$ has compact resolvent and generates a strongly continuous contraction semigroup on $\mathcal{X}$, by the Arendt-Batty-Lyubich-V\~u Theorem, see e.g.\ \cite[Theorem V.2.21]{EngelNagel_2000}, the semigroup is asymptotically stable if and only if
 \[
  \sigma_p(\mathcal{A}) \cap \ii \R
   = \emptyset.
 \]
Therefore, let $(\hat x, \beta) \in \dom(\mathcal{A}) \times \R$ be such that $\ii \beta \hat x = \mathcal{A} \hat x$.
Then, in particular
 \[
  0
   = \Re \sp{\ii \beta \hat x}{\hat x}_{\mathcal{X}}
   = \Re \sp{\mathcal{A} \hat x}{\hat x}_{\mathcal{X}}
   \leq - \kappa \abs{\mathfrak{R} x}_{\K^4}^2
   \leq 0,
 \]
i.e.\ $\mathfrak{R} x = 0$.
Then at least three of the four components of $((\H^1 x^1)(0), (\H^1 x^1)'(0))$ are zero.
Moreover, since the systems $\mathfrak{S} = (\mathfrak{A}_0, \mathfrak{B}, \mathfrak{C})$ and $\Sigma_c = (A_c, B_c, C_c, D_c)$ are both impedance passive and interconnected by standard feedback interconnection, it follows from Remark \ref{rem:on_assmpt_3.1} that
 \[
  0
   = \Re \sp{\ii \beta \hat x}{\hat x}_{\mathcal{X}}
   \leq - \sum_{j=0}^m \kappa_j \abs{D_c^j \mathfrak{C}^j x}_{U_c^j}^2,
 \]
which means that $D_c^j \mathfrak{C}^j x = 0$ for $j = 0, 1, \ldots, m$ as well.
Since $\ker D_c^j \subseteq \ker B_c^j$, this implies that $B_c^j \mathfrak{C}^j x = 0$ for $j = 0, 1, \ldots, m$ as well and, therefore,
 \[
  \ii \beta x_c^j
   = A_c^j x_c^j + B_c^j \mathfrak{C}^j x
   \quad \Rightarrow \quad
  x_c^j
   = (\ii \beta - A_c^j)^{-1} B_c^j \mathfrak{C}^j x
   = 0.
 \]
Then again, $\mathfrak{B}^j x = - C_c^j x_c^j - D_c^j \mathfrak{C}^j x = 0$ for $j = 0, 1, \ldots m$ follows at once.
By definition of $\mathfrak{B}^j$, this means that $(\H^j x^j)(1) = (\H^{j+1} x^{j+1})(0)$ and $(\H^j x^j)'(1) = (\H^{j+1} x^{j+1})'(0)$ for every junction $j = 1, \ldots, m-1$.

First, assume that $\beta \not= 0$.
After possible multiplication of $\hat x$ by some scalar $\alpha \in \C \setminus \{0\}$, we may and will assume that $(\H_1^1 x^1_1)(0) \geq 0$, $(\H^1_1 x^1_1)'(0) \geq 0$, $\ii \beta (\H^1_2 x^1_2)(0) \geq 0$ and $\ii \beta (\H^1_2 x^1_2)'(0) \geq 0$ (since three of these terms equal zero anyway).
By \cite[Lemma 4.2.9]{Augner_2016}, then either $\H^1 x^1 = 0$ on $[0,1]$, or $(\H_1^1 x^1_1)(\z) > 0$, $(\H^1_1 x^1_1)'(\z) > 0$, $\ii \beta (\H^1_2 x^1_2)(\z) > 0$ and $\ii \beta (\H^1_2 x^1_2)'(\z) > 0$ for all $\z \in (0,1]$.
Repeating this procedure and using the continuity conditions $(\H^j x^j)(1) = (\H^{j+1} x^{j+1})(0)$ and $(\H^j x^j)'(1) = (\H^{j+1} x^{j+1})'(0)$ then shows the following:
 \begin{align*}
  &\H^1 x^1 \equiv 0,
   \quad \Rightarrow \quad
   (\H^1 x^1)(1) = (\H^1 x^1)'(1) = 0
   \\
  &\quad \Rightarrow \quad
   (\H^2 x^2)(0) = (\H^2 x^2)'(0) = 0
   \quad \Rightarrow \quad
   \H^2 x^2 \equiv 0
   \\
   &\quad \Rightarrow \quad
   \H^j x^j \equiv 0, \, j = 1, \ldots, m
   \quad \Rightarrow \quad
   \hat x \equiv 0 \quad \text{ is not an eigenfunction},
   \\
   \intertext{and}
  &\H^1 x^1 \not\equiv 0
  \quad \Rightarrow \quad
  (\H^1_1 x^j_1)(1) > 0, \, (\H^1_1 x^j_1)'(1) > 0, \, \ii \beta (\H^1_2 x^1_2)(1) > 0, \, \ii \beta (\H^1_2 x^1_2)'(1) > 0
   \\
  &\quad \Rightarrow \quad
  (\H^2_1 x^2_1)(0) > 0, \, (\H^2_1 x^2_1)'(0) > 0, \, \ii \beta (\H^2_2 x^2_2)(0) > 0, \, \ii \beta (\H^2_2 x^2_2)'(0) > 0
   \\
  &\quad \Rightarrow \quad
  (\H^2_1 x^2_1)(\z) > 0, \, (\H^2_1 x^2_1)'(\z) > 0, \, \ii \beta (\H^2_2 x^2_2)(\z) > 0, \, \ii \beta (\H^2_2 x^2_2)'(\z) > 0
  \quad \text{for all } \z \in [0,1]
   \\
   &\quad \Rightarrow \quad
   (\H^j_1 x^j_1)(0) > 0, \, (\H^j_1 x^j_1)'(0) > 0, \, \ii \beta (\H^j_2 x^j_2)(0) > 0, \, \ii \beta (\H^j_2 x^j_2)'(0) > 0
   \quad \text{for } j = 1, \ldots, m
   \\
   &\quad \Rightarrow \quad
   (\H^j_1 x^j_1)(\z) > 0, \, (\H^j_1 x^j_1)'(\z) > 0, \, \ii \beta (\H^j_2 x^j_2)(\z) > 0, \, \ii \beta (\H^j_2 x^j_2)'(\z) > 0
   \quad \text{for all } \z \in [0,1], \, j = 1, \ldots, m
   \\
   &\quad \Rightarrow \quad
   \text{contradiction with the condition } \mathfrak{R} x = 0, \text{ hence, this case is impossible.}
 \end{align*}
This excludes the situation where $\sigma_p(A) \cap (\ii \R \setminus \{0\}) \neq \emptyset$.

For the case $\beta = 0$, one has $(\H^j x^j)'' = 0$ on $(0,1)$ for each $j = 1, \ldots, m$, but as above one sees that $(\H^{j+1} x^{j+1})(0) = (\H^j x^j)(1)$ and $(\H^{j+1} x^{j+1})'(0) = (\H^j x^j)'(1)$, so that
 \[
  ((\H^j x^j)(\z), (\H^j x^j)'(\z))
   = ((\H^1 x^1)(0), (\H^1 x^1)'(0))
   = ((\H^m x^m)(1), (\H^m x^m)'(1)),
   \quad
   \z \in [0,1], \, j = 1, \ldots, m,
 \]
but by the condition $\mathfrak{R} x = 0$ this can only be the case if $x = 0$, and then as before also $x_c = 0$,
Hence, we have shown that $\sigma(\mathcal{A}) \cap \ii \R = \emptyset$.
Asymptotic stability follows by the Arendt-Batty-Lyubich-V\~u-Theorem.\qed
\newline
In the language of \cite{AugnerJacob_2014} and as a by-product, the proof of Theorem \ref{thm:asymptotic_stability} shows the following:

\begin{lmm}
Let $\mathfrak{S}_0$, $\Sigma_c$ and $\mathfrak{R} = (\mathfrak{R}_1, \ldots, \mathfrak{R}_4)$ be as in Theorem \ref{thm:asymptotic_stability}.
 \begin{enumerate}
  \item
   The pair $(\mathfrak{A}_0, \mathfrak{R})$ has property \textrm{ASP}, i.e.\
 \[
  \forall \beta \in \R: \quad
  \ker(\mathfrak{A}_0 - \ii \beta) \cap \ker \mathfrak{R} = \{0\}.
 \]
 \item
  Let $\mathfrak{R}' = (\mathfrak{R}_1, \mathfrak{R}_2, \mathfrak{R}_3): \dom(\mathfrak{R}') = \dom(\mathfrak{R}) \rightarrow \K^3$, then
 \[
  \forall \beta \in \R \setminus \{0\}: \quad
  \ker(\mathfrak{A}_0 - \ii \beta) \cap \ker \mathfrak{R}' = \{0\}.
 \]
In particular, in the situation of Theorem \ref{thm:asymptotic_stability}, we have that $\sigma_p(\mathcal{A}) \cap \ii \R \subset \{0\}$, if the condition \eqref{eqn:diss_asymptoctic} is weakened to
 \[
  \Re \sp{\mathcal{A} (x,x_c)}{(x,x_c)}_{\mathcal{X}}
   \leq - \kappa \abs{\mathfrak{R}' x}_{\K^3}^2,
   \quad
   x \in \dom(\mathcal{A}).
 \]
 \end{enumerate}
\end{lmm}

Asymptotic stability may be seen as the first step towards exponential stability.
To even ensure uniform exponential stability, we adjust the setting by
 \begin{enumerate}
  \item
   imposing further restrictions on the boundary conditions at the left ($\z = \z_0 = 0$) and right end  ($\z = \z_m = L$) of the chain of beams and
  \item
   imposing monotonicity conditions on the parameter functions $\rho$ and $EI$ at the junction points $l^j \in (0,L)$, $j = 1, \ldots, m-1$, see condition \eqref{M} resp.\ condition \eqref{M'} below.
 \end{enumerate}

Condition \eqref{M} in the abstract port-Hamiltonian setting reads as

\begin{assumption}[Jump conditions in port-Hamiltonian formulation]
 We assume that
  \[
   \H^{j+1}(0) - \H^j(1)
   \quad \text{is positive semidefinite, for all }
   j = 1, \ldots, m-1.
   \tag{{\bf M'}}
   \label{M'}
  \]
\end{assumption}

\begin{thrm}[Uniform Exponential Stability]
 \label{thm:exp_stability}
  Let Assumption \ref{assmpt:passive_stable_controller} and conditions \eqref{R} and \eqref{M'} be satisfied and the operator $\mathcal{A}$ be such that
  \[
   \Re \sp{\mathcal{A} (x,x_c)}{(x,x_c)}_{\mathcal{X}}
    \leq - \kappa \abs{\mathfrak{R} x}_{\K^5}^2
  \]
 for some constant $\kappa > 0$ with $\mathfrak{R}: \dom(\mathfrak{A}) \rightarrow \K^5$ of the form
  \[
   \mathfrak{R} x
    = \left( \begin{array}{c} (\H^1 x^1)(0) \\ (\H^1_1 x^1_1)'(0) \, \text{or} \, (\H^1_2 x^1_2)'(0) \\ (\H^m_1 x^m_1)(1) \, \text{or} \, (\H^m_2 x^m_2)'(1) \\ (\H^m_1 x^m_1)'(1) \, \text{or} \, (\H^m_2 x^m_2)(1) \end{array} \right).
  \]
 Then, the $C_0$-semigroup $(\mathcal{T}(t))_{t \geq 0}$ generated by $\mathcal{A}$ is uniformly exponentially stable if and only if $(\mathcal{T}(t))_{t \geq 0}$ is asymptotically stable.\newline
 In particular, this is the case if $\sigma(A_c^j) \subseteq \C_0^-$ for $j = 0, 1, \ldots, m$ and
  \[
   \mathfrak{R} x
    = \left( \begin{array}{c} (\H^1 x^1)(0) \\ (\H^1_1 x^1_1)'(0) \\ (\H^m_1 x^m_1)(1) \, \text{or} \, (\H^m_2 x^m_2)'(1) \\ (\H^m_2 x^m_2)(1) \end{array} \right)
    \quad \text{or} \quad
   \mathfrak{R} x
    = \left( \begin{array}{c} (\H^1 x^1)(0) \\ (\H^1_2 x^1_2)'(0) \\ (\H^m_1 x^m_1)(1) \, \text{or} \, (\H^m_2 x^m_2)'(1) \\ (\H^m_1 x^m_1)'(1) \end{array} \right).
  \]
\end{thrm}

 \textbf{Proof.}
 First, we show that the pair $(\mathfrak{A}_0, (\mathfrak{R}, \mathfrak{B}))$ has property AIEP as introduced in Definition 2.8 of \cite{AugnerJacob_2014}, i.e.\ for every sequence $(x_n, \beta_n)_{n \geq 1} \subseteq \dom(\mathfrak{A}_0) \times \R$ with $\sup_{n \in \N} \norm{x_n}_X < \infty$ and $\abs{\beta_n} \rightarrow 0$ and such that
  \[
   \ii \beta_n x_n - \mathfrak{A}_0 x_n
    \rightarrow 0 \quad \text{in } X,
    \quad
   \mathfrak{R} x_n
    \rightarrow 0 \quad \text{in } \K^5,
    \quad
   \mathfrak{B} x_n
    \rightarrow 0 \quad \text{in } \K^{2m},
  \]
 it follows that $x_n \rightarrow 0$ in $X$.
 To this end, for $j = 1, \ldots, m$ fix functions $q^j \in C^2([0,1];\R)$ which will be specified below.
 By the proof of \cite[Proposition 4.3.19]{Augner_2016}, for every sequence $(x_n, \beta_n)_{n \geq 1} \subseteq \dom(\mathfrak{A}_0) \times \R$ with $\sup_{n \in \N} \norm{x_n}_X < \infty$, $\abs{\beta_n} \rightarrow 0$ and $\mathfrak{A}_0 x_n - \ii \beta_n x_n \rightarrow 0$ in $X$, it holds that
  \begin{align*}
   &\Re \sp{x^j_n}{(2 (q^j)' \H^j - q^j (\H^j)') x^j_n}_{L_2}
    \\
    &= o(1)
     - 2 \Re \left[ \sp{- (\H^j_2 x^j_{n,2})'(\z)}{\frac{\ii q^j(\z)}{\beta_n} (\H^j_1 x^j_{n,1})'(\z)}_{\K} \right]_0^1
     + \left[ \sp{x^j_n(\z)}{q^j \H^j(\z) x^j_n(\z)}_{\K} \right]_0^1
     \\
     &\quad
     + \left[ \Re \sp{- (\H^j_2 x^j_{n,2})'(\z)}{\frac{\ii (q^j)'(\z)}{\beta_n} (\H^j_1 x^j_{n,1})(\z)}_{\K} \right]_0^1 
     - \left[ \Re \sp{-(\H^j_2 x^j_{n,2})(\z)}{\frac{\ii (q^j)'}{\beta_n} (\H^j_1 x^j_{n,1})'(\z)}_{\K} \right]_0^1
  \end{align*}
where $o(1)$ denotes further terms which tend to zero as $n \rightarrow \infty$.
(Note that there is a typo in \cite[equation (4.27)]{Augner_2016}: There actually should be a minus sign in front of the last line of the equation.)
Summing up these equalities and writing $Q(\z) = \operatorname{diag}\, (q^j(\z))_{j=1}^m$, we find that
{\allowdisplaybreaks[1]
  \begin{align*}
   &\Re \sp{x_n}{(2 Q' \H - Q \H') x_n}_{L_2}
    \\
    &= o(1)
     - 2 \Re \left[ \sp{- (\H^m_2 x^m_{n,2})'(1)}{\frac{\ii q^m(1)}{\beta_n} (\H^m_1 x^m_{n,1})'(1)}_{\K} \right]
     + \left[ \sp{x^m_n(1)}{q^m \H^m(1) x^m_n(1)}_{\K} \right]_0^1
     \\
     &\quad
     +  \Re \sp{- (\H^m_2 x^m_{n,2})'(1)}{\frac{\ii (q^m)'(1)}{\beta_n} (\H^m_1 x^m_{n,1})(1)}_{\K}
     - \Re \sp{-(\H^m_2 x^m_{n,2})(1)}{\frac{\ii (q^m)'(1)}{\beta_n} (\H^m_1 x^m_{n,1})'(1)}_{\K}
     \\
     &\quad
     + 2 \Re \sp{- (\H^1_2 x^1_{n,2})'(0)}{\frac{\ii q^1(0)}{\beta_n} (\H^1_1 x^1_{n,1})'(0)}_{\K} 
     - \sp{x^1_n(0)}{q^1 \H^1(0) x^1_n(0)}_{\K}
     \\
     &\quad
     - \Re \sp{- (\H^1_2 x^1_{n,2})'(0)}{\frac{\ii (q^1)'(0)}{\beta_n} (\H^1_1 x^1_{n,1})(0)}_{\K}
     + \Re \sp{-(\H^1_2 x^1_{n,2})(0)}{\frac{\ii (q^1)'(0)}{\beta_n} (\H^1_1 x^1_{n,1})'(0)}_{\K}
     \\
     &\quad
     - 2 \sum_{j=1}^{m-1} \Re \left[ \sp{- (\H^j_2 x^j_{n,2})'(1)}{\frac{\ii q^j(1)}{\beta_n} (\H^j_1 x^j_{n,1})'(1)}_{\K} - \sp{- (\H^{j+1}_2 x^{j+1}_{n,2})'(0)}{\frac{\ii q^{j+1}(0)}{\beta_n} (\H^{j+1}_1 x^{j+1}_{n,1})'(0)}_{\K} \right]
     \\
     &\quad
     +  2 \sum_{j=1}^{m-1} \left[ \sp{x^j_n(1)}{q^j \H^j(1) x^j_n(1)}_{\K} - \sp{x^{j+1}_n(0)}{q^{j+1} \H^{j+1}(0) x^{j+1}_n(0)}_{\K} \right]
     \\
     &\quad
     -  2 \sum_{j=1}^{m-1} \left[ \Re \sp{(\H^j_2 x^j_{n,2})'(1)}{\frac{\ii (q^j)'(1)}{\beta_n} (\H^j_1 x^j_{n,1})(1)}_{\K} \hspace{-0.15cm} - \Re \sp{(\H^{j+1}_2 x^{j+1}_{n,2})'(0)}{\frac{\ii (q^{j+1})'(0)}{\beta_n} (\H^{j+1}_1 x^{j+1}_{n,1})(0)}_{\K} \right]
     \\
     &\quad
     +  2 \sum_{j=1}^{m-1} \left[ \Re \sp{(\H^j_2 x^j_{n,2})(1)}{\frac{\ii (q^j)'(1)}{\beta_n} (\H^j_1 x^j_{n,1})'(1)}_{\K} \hspace{-0.15cm} - \Re \sp{(\H^{j+1}_2 x^{j+1}_{n,2})(0)}{\frac{\ii (q^{j+1})'(0)}{\beta_n} (\H^{j+1}_1 x^{j+1}_{n,1})'(0)}_{\K} \right]_0^1.
  \end{align*}
 }We show that for a suitable choice of the functions $q^j$ and under the imposed dissipation and monotonicity assumptions, all the terms on the right hand side vanish as $n \rightarrow \infty$.
 To this end, observe that $\ii \beta_n x_n - \mathfrak{A} x_n \rightarrow 0$ in $X$ and $(x_n)_{n \geq 1} \subseteq X$ is bounded whereas $\abs{\beta_n} \rightarrow \infty$, so that $\frac{\mathfrak{A}_0 x_n}{\beta_n} = P_2 (\H x_n)'' \rightarrow 0$ in $X$, thus $(\H x_n)'' \rightarrow 0$ in $X$.
 From \cite[Lemma 2.15]{AugnerJacob_2014}, an interpolation argument and embedding theorems for Sobolev Slobodetskii spaces into $C^k$-spaces we conclude that
  \begin{equation}
   \frac{\H x_n}{\beta_n} \rightarrow 0,
    \quad
    \text{in } C^1([0,1];\K^{2m}).
    \label{eqn:C^1-o(1)}
  \end{equation}
 Only from here on, we additionally assume that $\mathfrak{R} x_n \rightarrow 0$.
 \begin{enumerate}
  \item
   The terms
    \[
     - 2 \Re \left[ \sp{- (\H^m_2 x^m_{n,2})'(1)}{\frac{\ii q^m(1)}{\beta_n} (\H^m_1 x^m_{n,1})'(1)}_{\K} \right]
     + \left[ \sp{x^m_n(1)}{q^m \H^m(1) x^m_n(1)}_{\K} \right]_0^1
    \]
   not only tend to zero, but equal zero for all $n \in \N$, if we choose the function $q^m$ such that
    \begin{equation}
     q^m(1) = 0.
     \label{cond:q-1}
    \end{equation}
   \item
    For the third term, we have
     \begin{align*}
      \abs{ \Re \sp{- (\H^m_2 x^m_{n,2})'(1)}{\frac{\ii (q^m)'(1)}{\beta_n} (\H^m_1 x^m_{n,1})(1)}_{\K} }
       &\leq |(q^m)'(1)| \frac{1}{\abs{\beta_n}} \abs{(\H^m_2 x^m_{n,2})'(1))} \abs{(\H^m_1 x^m_{n,1})}
       \rightarrow 0
     \end{align*}
    as from $\mathfrak{R} x_n \rightarrow 0$ and the definition of $\mathfrak{R}$, at least one of the terms $\abs{(\H^m_2 x^m_{n,2})'(1))}$ or $\abs{(\H^m_1 x^m_{n,1})(1)}$ tends to zero and by \eqref{eqn:C^1-o(1)} in any case $\frac{1}{\abs{\beta_n}} \abs{(\H^m_1 x^m_{n,1})(1)}$ and $\frac{1}{\abs{\beta_n}} \abs{(\H^m_2 x^m_{n,2})'(1))}$ tend to zero.
   \item
    Similarly,
     \begin{align*}
      \abs{ \Re \sp{-(\H^m_2 x^m_{n,2})(1)}{\frac{\ii (q^m)'(1)}{\beta_n} (\H^m_1 x^m_{n,1})'(1)}_{\K} }
       &\leq |(q^m)'(1)| \frac{1}{\abs{\beta_n}} \abs{(\H^m_2 x^m_{n,2})(1))} \abs{(\H^m_1 x^m_{n,1})'}
       \rightarrow 0
     \end{align*}
    as $\abs{(\H^m_2 x^m_{n,2})(1))}$ or $\abs{(\H^m_1 x^m_{n,1})'(1)}$ tends to zero due to $\mathfrak{R} x_n \rightarrow 0$ and the definition of $\mathfrak{R}$.
   \item
   Since at least one of the terms $(\H^1_2 x^1_{n,2})'(0)$ or $(\H^1_1 x^1_{n,1})'(0)$ tends to zero as well, by the same reasoning
    \[
     \abs{ \Re \sp{- (\H^1_2 x^1_{n,2})'(0)}{\frac{\ii q^1(0)}{\beta_n} (\H^1_1 x^1_{n,1})'(0)}_{\K} }
      \rightarrow 0.
    \]
  \item
   For the terms
    \begin{align*}
     - \left[ \sp{x^1_n(0)}{q^1 \H^1(0) x^1_n(0)}_{\K} \right]
     - \Re \sp{- (\H^1_2 x^1_{n,2})'(0)}{\frac{\ii (q^1)'(0)}{\beta_n} (\H^1_1 x^1_{n,1})(0)}_{\K}
     \\
     - \Re \sp{-(\H^1_2 x^1_{n,2})(0)}{\frac{\ii (q^1)'(0)}{\beta_n} (\H^1_1 x^1_{n,1})'(0)}_{\K}
     &\rightarrow 0
    \end{align*}
   we use that $(\H^1 x^1_n)(0) \rightarrow 0$ since $\mathfrak{R} x_n \rightarrow 0$ and then $(\H^1 x^1_n)(0) \rightarrow 0$.\newline
   This concludes the discussion of the boundary terms at the left and right end of the chain of beams.
   For the junction points, we use the continuity and jump conditions, and obtain the following.
  \item
   For each $j = 1, \ldots, m-1$ we have that
    \begin{align*}
     &+ \Re \sp{- (\H^j_2 x^j_{n,2})'(1)}{\frac{\ii (q^j)'(1)}{\beta_n} (\H^j_1 x^j_{n,1})(1)}_{\K} 
      - \Re \sp{- (\H^{j+1}_2 x^{j+1}_{n,2})'(0)}{\frac{\ii (q^{j+1})'(0)}{\beta_n} (\H^{j+1}_1 x^{j+1}_{n,1})(0)}_{\K}
      \\
      &- \Re \sp{- (\H^j_2 x^j_{n,2})'(1)}{\frac{\ii (q^j)'(1)}{\beta_n} (\H^j_1 x^j_{n,1})(1)}_{\K} + \Re \sp{- (\H^{j+1}_2 x^{j+1}_{n,2})'(0)}{\frac{\ii (q^{j+1})'(0)}{\beta_n} (\H^{j+1}_1 x^{j+1}_{n,1})(0)}_{\K}
      \\
      &=
       \Re \sp{ \mathfrak{B}^j x_n}{ \frac{\ii (q^{j+1})'(0)}{\beta_n} \mathfrak{C}^j x_n}
    \end{align*}
   if we choose the functions $q^j$ and $q^{j+1}$ such that their derivatives at the junction points match, i.e.\
    \begin{equation}
     (q^j)'(1)
      = (q^{j+1})'(0),
      \quad
      j = 1, \ldots, m-1.
      \label{cond:q-2}
    \end{equation}
   One may easily check this claim for all the cases of the control and observation maps $\mathfrak{B}^j$ and $\mathfrak{C}^j$ we allowed for.
    Moreover, since $\mathfrak{B} x_n \rightarrow 0$ by assumption and $\frac{\mathfrak{C} x_n}{\beta_n} \rightarrow 0$ by \eqref{eqn:C^1-o(1)}, we find that these terms tend to zero as $n \rightarrow \infty$ as well:
      \[
       \Re \sp{ \mathfrak{B}^j x}{ \frac{\ii (q^{j+1})'(0)}{\beta_n} \mathfrak{C}^j x}
        \rightarrow 0.
      \]
   \item
    To handle the terms
     \begin{align*}
      &\Re \sp{- (\H^j_2 x^j_{n,2})'(1)}{\frac{\ii q^j(1)}{\beta_n} (\H^j_1 x^j_{n,1})'(1)}_{\K} 
       - \Re \sp{- (\H^{j+1}_2 x^{j+1}_{n,2})'(0)}{\frac{\ii q^{j+1}(0)}{\beta_n} (\H^{j+1}_1 x^{j+1}_{n,1})'(0)}_{\K}
       \\
      &= \Re \sp{- (\H^j_2 x^j_{n,2})'(1) + (\H^{j+1}_2 x^{j+1}_{n,2})'(0)}{\frac{\ii q^{j+1}(0)}{\beta_n} (\H^{j+1}_1 x^{j+1}_{n,1})'(0)}_{\K}
       \\
       &\quad
       - \Re \sp{- (\H^{j+1}_2 x^{j+1}_{n,2})'(0)}{\frac{\ii q^{j+1}(0)}{\beta_n} \left( (\H^{j+1}_1 x^{j+1}_{n,1})'(0) - (\H^j_1 x^j_{n,1})'(1) \right)}_{\K}
       \\
      &= \Re \sp{- (\H^j_1 x^j_{n,1})'(1) + (\H^{j+1}_1 x^{j+1}_{n,1})'(0)}{\frac{\ii q^{j+1}(0)}{\beta_n} (\H^{j+1}_2 x^{j+1}_{n,2})'(0)}_{\K}
       \\
       &\quad
       - \Re \sp{- (\H^{j+1}_1 x^{j+1}_{n,1})'(0)}{\frac{\ii q^{j+1}(0)}{\beta_n} \left( (\H^{j+1}_2 x^{j+1}_{n,2})'(0) - (\H^j_1 x^j_{n,2})'(1) \right)}_{\K},
     \end{align*}
    we note that both terms $(\H^j_1 x^j_1)'(1) - (\H_1^{j+1} x^{j+1}_1)'(0)$ and $(\H^j_2 x^j_2)'(1) - (\H_2^{j+1} x^{j+1}_2)'(0)$ tend to zero (if the corresponding term is a component of $\mathfrak{B}^j x$) or even equal zero (if the corresponding term constitutes a component of $\mathfrak{B}_0^j x$).
    As $\frac{1}{\beta} (\H_1^{j+1} x^{j+1}_1)'(0)$ and $\frac{1}{\beta_n} (\H_2^{j+1} x_2^{j+1})'(0)$ tend to zero by \eqref{eqn:C^1-o(1)}, it follows that
     \begin{align*}
      &\Re \sp{- (\H^j_2 x^j_{n,2})'(1)}{\frac{\ii q^j(1)}{\beta_n} (\H^j_1 x^j_{n,1})'(1)}_{\K} 
       \\
       &- \Re \sp{- (\H^{j+1}_2 x^{j+1}_{n,2})'(0)}{\frac{\ii q^{j+1}(0)}{\beta_n} (\H^{j+1}_1 x^{j+1}_{n,1})'(0)}_{\K}
       \longrightarrow 0,
      \quad
      n \rightarrow \infty.
     \end{align*}
  \item
   Lastly, the terms
    \begin{align*}
     &\sp{x^j_n(1)}{(q^j \H^j)(1) x^j_n(1)}_{\K} - \sp{x^{j+1}_n(0)}{(q^{j+1} \H^{j+1})(0) x^{j+1}_n(0)}_{\K}
      \\
      &= q^{j+1}(0) \sp{x^j_n(1)}{\H^j(1) x^j_n(1)}_{\K} - \sp{x^{j+1}_n(0)}{\H^{j+1}(0) x^{j+1}_n(0)}_{\K}
    \end{align*}
   have to be handled, where due to previous restrictions on the choice of the functions $q^j$ and $q^{j+1}$, necessarily $q^j(1) = q^{j+1}(0)$.
   For the moment, we leave these terms as they are and employ the monotonicity conditions on $\H^j(1) - \H^{j+1}(0)$ later on.
 \end{enumerate}
Putting things together, and choosing the functions $q^j$ such that $q^j(\z) = q(j-1+\z)$ for some function $q \in C^2([0,m];\R)$ with $q(m) = 0$, we find that
 \begin{align*}
  &\Re \sp{x_n}{(2 Q' \H - Q \H') x_n}_{L_2}
    \\
    &= 2 \sum_{j=1}^{m-1} q(j) \left[ \sp{x^j_n(1)}{\H^j(1) x^j_n(1)}_{\K} - \sp{x^{j+1}_n(0)}{\H^{j+1}(0) x^{j+1}_n(0)}_{\K} \right]
    + o(1).
 \end{align*}
We will now choose the function $q$ such that the term on the left hand side defines an equivalent inner product on $L_2(0,1;\K^{2m})$ and to this end choose $q \in C^2([0,m];\R)$ with $q(m) = 0$ and $q' > 0$, $q \leq 0$ such that
 \[
  2 (q^j)' \H^j - q^j (\H^j)'
   \geq 2 m_0 (q^j)' + q^j M_1
   \geq \varepsilon_0
 \]
where
 \begin{align*}
  m_0
   &:= \sup \{ \varepsilon > 0: \H^j(\z) - \varepsilon I \text{ positive semidefinite for a.e.\, } \z \in (0,1), j = 1, \ldots, m \},
   \\
  M_1
   &:= \inf \{ \varepsilon > 0: \varepsilon I - (\H^j)'(\z) \text{ positive semidefinite for a.e.\, } \z \in (0,1), j = 1, \ldots, m \}
 \end{align*}
and $\varepsilon_0 > 0$.
We then have the estimate
 \[
  \varepsilon_0 \norm{x_n}_{L_2}^2
   \leq 2 \sum_{j=1}^{m-1} q(j) \left[ \sp{x^j_n(1)}{\H^j(1) x^j_n(1)}_{\K} - \sp{x^{j+1}_n(0)}{\H^{j+1}(0) x^{j+1}_n(0)}_{\K} \right]
    + o(1).
 \]
Since for every $j = 1, \ldots, m-1$, the term $q(j) < 0$ is strictly negative, the right hand side is less or equal $o(1)$, if for each junction $j = 1, \ldots, m-1$, we have that
 \begin{align*}
  &\sp{x^j_n(1)}{\H^j(1) x^j_n(1)}_{\K}
   - \sp{x^{j+1}_n(0)}{\H^{j+1}(0) x^{j+1}_n(0)}_{\K}
   \\
   &= \sp{(\H^j x^j_n)(1) - (\H^{j+1} x^{j+1}_n)(0)}{(\H^j (1)^{-1} - \H^{j+1}(0)^{-1}) \left(  (\H^j x^j_n)(1) - (\H^{j+1} x^{j+1}_n)(0) \right)}_{\K}
    \\
    &\quad
    + 2 \sp{(\H^{j+1} x^{j+1}_n)(0)}{(\H^j(1)^{-1} - \H^{j+1}(0)^{-1}) (\H^j x^j_n)(1)} 
    \\
   &= \frac{3}{2} \sp{(\H^j x^j_n)(1) - (\H^{j+1} x^{j+1}_n)(0)}{(\H^j (1)^{-1} - \H^{j+1}(0)^{-1}) \left(  (\H^j x^j_n)(1) - (\H^{j+1} x^{j+1}_n)(0) \right)}_{\K}
    \\
    &\quad
    + \frac{1}{2} \sp{\left(  (\H^j x^j_n)(1) + (\H^{j+1} x^{j+1}_n)(0) \right)}{(\H^j(1)^{-1} - \H^{j+1}(0)^{-1}) \left(  (\H^j x^j_n)(1) + (\H^{j+1} x^{j+1}_n)(0) \right)}
 \end{align*}
 is a sum of $o(1)$-terms.
 By assumption, $\mathfrak{B}^j x_n \rightarrow 0$ and $\mathfrak{B}_0^j x_n = 0$.
 By this, the term $(\H^j x^j_n)(1) - (\H^{j+1} x^{j+1}_n)(0)$ and, hence, the first term of the latter sum converge to zero as $n \rightarrow \infty$.
 For the second term, we can at least ensure that it is non-negative, by using the structural monotonicity assumption that
  \[
   \H^j(1)^{-1} - \H^{j+1}(0)^{-1}
    \quad
    \text{is positive semidefinite}.
  \]
 (Recall that in terms of $\rho$ and $EI$ this means that $\rho(l^j-) \geq \rho(l^j+)$ and $(EI)(l^j-) \leq (EI)(l^j+)$.)
 Then,
  \[
   \sp{x_n}{(2 Q' \H - Q \H') x_n}_{L_2}
    \rightarrow 0
  \]
 and since this inner product induces a norm which is equivalent to the standard $L_2$-norm $\norm{\cdot}_{L_2}$ and the energy norm $\norm{\cdot}_X$ on $X$, this implies that $x_n \rightarrow 0$ in $X$. We have successfully proved property AIEP.
 Let us state this preliminary result as
 
 \begin{prpstn}
 \label{prop:AIEP}
  Let the conditions of Theorem \ref{thm:exp_stability} be satisfied.
  Then, the pair $(\mathfrak{A}_0, (\mathfrak{R}, \mathfrak{B}))$ has property AIEP, i.e.\ for every sequence $(x_n, \beta_n)_{n \geq 1} \subseteq \dom(\mathfrak{A}_0) \times \R$ with
   \begin{align*}
    \sup_{n \in \N} \norm{x_n}_X 
     &< \infty,
     \\
    \abs{\beta_n}
     &\rightarrow 0,
     \\
    (\mathfrak{A}_0 - \ii \beta_n) x_n
     &\rightarrow 0,
     \\
    \mathfrak{R} x_n, \mathfrak{B} x_n
     &\rightarrow 0
   \end{align*}
  it follows that
   \[
    \norm{x_n}_{L_2}
     \rightarrow 0.
   \]
 \end{prpstn} 

 \textbf{Proof of Theorem \ref{thm:exp_stability} (continued).}
 We are now in position to show the assertion of Theorem \ref{thm:exp_stability}.
 Assume that $(\mathcal{T}(t))_{t \geq 0}$ is asymptotically stable, so that $\sigma(\mathcal{A}) = \sigma_p(\mathcal{A}) \subseteq \R$ by the Arendt-Batty-Lyubich-V\~u Theorem.
 Thanks to the Gearhart-Pr\"uss-Huang Theorem, see e.g.\ \cite[Theorem V.1.11]{EngelNagel_2000}, it suffices to prove that
  \[
   \sup_{\beta \in \R} \norm{(\ii \beta - \mathcal{A})^{-1}}_{\mathcal{X}}
    < \infty
  \]
 and this is equivalent to the statement
  \[
   \left.
    \begin{array}{c}
     (\hat x_n, \beta_n)_{n \geq 1} \subseteq \dom(\mathcal{A}) \times \R
     \\
     \sup_{n \in \N} \norm{\hat x_n}_{\mathcal{X}} < \infty,
     \quad
     \abs{\beta_n} \rightarrow \infty,
     \quad
     \mathcal{A} \hat x_n - \ii \beta_n \hat x_n
      \rightarrow 0
     \end{array}
   \right\}
    \quad \Rightarrow \quad
    \norm{\hat x_n}_{\mathcal{X}} \rightarrow 0,
  \]
 where we write $\hat x_n = (x_n, x_{n,c}) \in X \times X_c$; see, e.g.\ \cite[Remark 2.7]{AugnerJacob_2014}.
 Let $(\hat x_n, \beta_n)_{n \geq 1}$ be such a sequence. Then, in particular
  \[
   0
    \leftarrow  \Re \sp{(\mathcal{A} - \ii \beta_n) \hat x_n}{\hat x_n}
    = \Re \sp{\mathcal{A} \hat x_n}{\hat x_n}
    \leq - \kappa \abs{\mathfrak{R} x_n}_{\K^5}
  \]
 and hence $\mathfrak{R} x_n \rightarrow 0$.
 Having property AIEP for the pair $(\mathfrak{A}_0, (\mathfrak{R}, \mathfrak{B}))$ at hand, we show that
  \begin{enumerate}
   \item
    $x_{c,n} \rightarrow 0$ in $X_c$ and
   \item
    $\mathfrak{B} x_n \rightarrow 0$ in $\K^{2(m+1)}$
  \end{enumerate}
 By property AIEP of the pair $(\mathfrak{A}_0, (\mathfrak{R}, \mathfrak{B}))$, we may then conclude $x_n \rightarrow 0$ in $X$ and hence $\hat x_n \rightarrow 0$ in $\mathcal{X}$.
 
 By \eqref{eqn:C^1-o(1)}, the term $\frac{\mathfrak{C} x_n}{\beta_n}$ tends to zero as $n \rightarrow \infty$ and, hence, we obtain that
  \[
   x_{c,n}
    = (\ii \beta_n - A_c)^{-1} (\mathfrak{C} x_n + o(1))
    = \beta_n (\ii \beta_n - A_c)^{-1}  \left( \frac{\mathfrak{C} x_n}{\beta_n} + \frac{1}{\abs{\beta_n}} o(1) \right)
    \rightarrow 0,
  \]
 as the resolvent operators $\norm{\beta (\ii \beta - A_c)^{-1}} = \norm{ (\ii - A_c/\beta)^{-1}}$ are uniformly bounded for $\beta \in \R$.
 Moreover, by Assumption \ref{assmpt:passive_stable_controller} and Remark \ref{rem:on_assmpt_3.1}
  \[
   0
    \leftarrow \Re \sp{(\mathcal{A} - \ii \beta_n) \hat x_n}{\hat x_n}_{\mathcal{X}}
    \leq - \sum_{j=0}^m \kappa_j \abs{D_c^j \mathfrak{C}^j x_n}^2,
  \]
 so that $D_c^j \mathfrak{C}^j x_n \rightarrow 0$.
 Since $\ker D_c^j \subseteq \ker B_c^j$ this implies that $B_c^j \mathfrak{C}^j x_n \rightarrow 0$ as well, thus
  \[
   \mathfrak{B} x_n
    = - (C_c (\ii \beta_n - A_c)^{-1} B_c + D_c) \mathfrak{C} x_n + \frac{o(1)}{\abs{\beta_n}}
    \rightarrow 0.
  \]
 By property AIEP of the pair $(\mathfrak{A}_0, (\mathfrak{R}, \mathfrak{B}))$ as stated in Proposition \ref{prop:AIEP}, it follows that $x_n \rightarrow 0$ in $X$, therefore, $\hat x_n \rightarrow 0$ in $\mathcal{X}$ and this concludes the proof.\qed

Let us translate the conditions on $\mathfrak{R}$ into suitable (for exponential energy decay) choices of the conservative boundary conditions at the right end of the chain.

 \begin{rmrk}
  Our approach covers all the conservative boundary conditions mentioned in Subsection 4.1 of \cite{ChenEtAl_1987}:
   \begin{enumerate}
    \item
     $\omega(t,1) = (EI \omega_{\z\z})(t,1)$ (simply supported or pinned right end),
    \item
     $\omega_{\z\z}(t,1) = (EI \omega_{\z\z})_\z(t,1) = 0$ (free right end),
    \item
     $\omega_\z(t,1) = (EI \omega_{\z\z})_\z(t,0) = 0$ (shear hinge right end),
    \item
     $\omega_t(t,1) = \omega_\z(t,1) = 0$ (clamped left end),
    \item
     $\omega_t(t,1) = (EI \omega_{\z\z})(t,1) = 0$,
    \item
     $\omega_{t\z}(t,1) = (EI \omega_{\z\z})_\z(t,1) = 0$.
   \end{enumerate}
  In the energy state space formulation we used for the proof of well-posedness, asymptotic and exponential stability these are actually only four cases:
  In the energy state space formulation one does neither distinct between the cases $\omega(t,1) = 0$ and $\omega_t(t,1) = 0$ (i.e.\ $\omega(t,1) = c$), cf.\ the first and fifth case in the list just above, nor between the cases $\omega_\z(t,1) = 0$ and $\omega_{t\z}(t,1) = 0$ (i.e.\ $\omega_\z(t,1) = c$), cf.\ the third and last case.
  In energy state space the conditions above therefore read as:
   \begin{enumerate}
    \item
     $(\H^m x^m)(1) = 0$,
    \item
     $(\H^m_2 x^m_2)(1) = (\H^m_2 x^m_2)'(1) = 0$,
    \item
     $(\H^m x^m)'(1) = 0$,
    \item 
     $(\H^m_1 x^m_1)(1) = (\H^m_1 x^m_1)'(1) = 0$,
   \end{enumerate}
 \end{rmrk} 
 
To make the formulation of the stability results more digestible, we introduce the following

\begin{assumption}[Boundary and interconnection conditions]
At the left and right end and at the junction points we assume the following:
\newline
Condition {\bf (D)}:
  \begin{align*}
  \left( \begin{array}{c} (EI \omega_{\z\z})(0) \\ - (EI \omega_{\z\z})_\z(0) \end{array} \right)
   &= - K_0 \left( \begin{array}{c} \omega_{t\z}(t,\z) \\ \omega_t(t,\z) \end{array} \right)
   \label{eqn:dissipative_end}
   \\
  \intertext{for some $K_0 \in \K^{2 \times 2}$ with}
  K_0
   &= \left( \begin{array}{cc} k^0_{11} & 0 \\ 0 & 0 \end{array} \right)
   \text{ for some } k^0_{11} > 0
   \quad
   \text{or}
   \quad
   \Sym(K_0) \quad \text{is symmetric positive definite},
  \end{align*}
 Condition {\bf (C)}:
     \begin{enumerate}
    \item
     $\omega_t(t,1) = (EI \omega_{\z\z})(t,1) = 0$ (simply supported or pinned right end) and $\Sym (K_0)$ is positive definite, or
    \item
     $\omega_{\z\z}(t,1) = (EI \omega_{\z\z})_\z(t,1) = 0$ (free right end) and $\Sym(K_0)$ is positive definite, or
    \item
     $\omega_{t\z}(t,1) = (EI \omega_{\z\z})_\z(t,0) = 0$ (shear hinge right end), or
    \item
     $\omega_t(t,1) = \omega_{t\z}(t,1) = 0$.
   \end{enumerate}
  Condition {\bf (I)}: For $j = 1, \ldots, m-1$,
 \begin{enumerate}
  \item
   $\omega_t(t,l^j-)  = \omega_t(t,l^j +)$, $\omega_{t\z}(t,l^j-)  = \omega_{t\z}(t,l^j +)$ and
   \[
    \left( \begin{array}{c} - (EI \omega_{\z\z})_\z(t,l^j-) + (EI \omega_{\z\z})_\z(t,l^j+) \\ (EI \omega_{\z\z})(t,l^j-) - (EI \omega_{\z\z})(t,l^j+) \end{array} \right) = - K^j \left( \begin{array}{c} \omega_t(t,l^j) \\ \omega_{t\z}(t,l^j) \end{array} \right)
   \]
  \item
   $\omega_t(t,l^j-)  = \omega_t(t,l^j +)$, $(EI \omega_{\z\z})(t,l^j-) = (EI \omega_{\z\z})(t,l^j+)$ and
   \[
    \left( \begin{array}{c} - (EI \omega_{\z\z})_\z(t,l^j-) + (EI \omega_{\z\z})_\z(t,l^j+) \\ \omega_{t\z}(t,l^j-) - \omega_{t\z}(t, l^j+) \end{array} \right) = - K^j \left( \begin{array}{c} \omega_t(t,l^j) \\ (EI \omega_{\z\z})(t,l^j) \end{array} \right)
   \]
  \item
   $- (EI \omega_{\z\z})_\z(t,l^j-) = - (EI \omega_{\z\z})_\z(t,l^j+)$, $\omega_{t\z}(t,l^j-)  = \omega_{t\z}(t,l^j +)$ and
   \[
    \left( \begin{array}{c} \omega_t(t,l^j-) - \omega_t(t,l^j+) \\ (EI \omega_{\z\z})(t,l^j-) - (EI \omega_{\z\z})(t,l^j+) \end{array} \right) = - K^j \left( \begin{array}{c} - (EI \omega_{\z\z})_\z(t,l^j) \\ \omega_{t\z}(t,l^j) \end{array} \right)
   \]
  \item
   $- (EI \omega_{\z\z})_\z(t,l^j-) = - (EI \omega_{\z\z})_\z(t,l^j+)$, $(EI \omega_{\z\z})(t,l^j-) = (EI \omega_{\z\z})(t,l^j+)$ and
   \[
    \left( \begin{array}{c} \omega_t(t,l^j-) - \omega_t(t,l^j+) \\ \omega_{t\z}(t,l^j-) - \omega_{t\z}(t, l^j+) \end{array} \right) = - K^j \left( \begin{array}{c} - (EI \omega_{\z\z})_\z(t,l^j) \\ (EI \omega_{\z\z})(t,l^j) \end{array} \right)
   \]
 \end{enumerate}
where $K^j \in \K^{2 \times 2}$ is a diagonal, positive semidefinite, symmetric matrix or has a positive definite symmetric part $\Sym(K^j) > 0$.
\end{assumption}

The exponential stability result for static boundary feedback and interconnection conditions, then translates as

 \begin{thrm}[Stability of Serially Connected Euler-Bernoulli Beams]
 \label{exa:EB}
 Let conditions \eqref{R}, \eqref{M}, {\bf (C)}, {\bf (D)} and {\bf (I)} be satisfied.
Then, for every initial datum
 \[
  (\omega(0,\cdot), \omega_t(0,\cdot))
   = (\omega_0, \omega_1)
   \in H^2((0,L) \setminus \{l^j\}_{j=1}^m )^2,
 \]
there is a unique strong solution
 \[
  \omega \in C([0,\infty); H^2((0,L) \setminus \{l^j\}_{j=1}^m)
 \]
of the Euler-Bernoulli-Beam system \eqref{EB} with the imposed boundary and interconnection conditions in {\bf (C)}, {\bf (D)} and {\bf (I)}.
 The solution depends continuously on the initial data $(\omega_0, \omega_1)$, and there are constants $M \geq 1$ and $\eta < 0$, independent of the initial data, such that the energy
 \[
  H(t)
   := \int_0^t \rho(\z) \abs{\omega_t(t,\z)}^2 + EI(\z) \abs{\omega_{\z\z}(t,\z)}^2 \dd \z
 \]
decays uniformly exponentially, i.e.\
 \[
  H(t)
   \leq M \ee^{\eta t} H(0),
   \quad
   t \geq 0.
 \]
 \end{thrm}
 
 \textbf{Proof.}
 In either case, we have
  \[
   \Re \sp{\mathcal{A} \hat x}{\hat x}
    \leq - \kappa \left( \abs{(\H^1_2 x^1_2)'(0)}^2 + \abs{(\H^1 x^1)'(0)}^2 \right),
    \quad
    x \in \dom(\mathcal{A})
  \]
 for some $\kappa > 0$.
 For the first and fourth case, asymptotic stability follows by Theorem \ref{thm:asymptotic_stability}, whereas for the second and third case the case that $0 \in \sigma_p(\mathcal{A})$ has to be excluded by demanding that $\Sym(K_0) > 0$, so that
  \[
   \Re \sp{\mathcal{A} \hat x}{\hat x}
    \leq - \kappa \left( \abs{(\H^1 x^1)'(0)}^2 + \abs{(\H^1 x^1)'(0)}^2 \right),
    \quad
    x \in \dom(\mathcal{A}).
  \]
 In all cases, uniform exponential stability follows by Theorem \ref{thm:exp_stability}.\qed
 
 \begin{rmrk}
  If, in the second and third case of Theorem \ref{exa:EB}, one considers the dissipative feedback with $K_0 = \operatorname{diag}(k^0_{11}, 0)$ for some $k^0_{11} > 0$ at the left end, either $\sigma_p(\mathcal{A}) \cap \ii \R = \emptyset$, e.g.\ by suitable damping in one of the junction points $l^j$, and then by Theorem \ref{thm:exp_stability} the system is again uniformly exponentially stable, or $\sigma_p(\mathcal{A}) \cap \ii \R = \{0\}$. (This follows from the proof of Theorem \ref{thm:asymptotic_stability} which shows that $\sigma_p(\mathcal{A}) \cap \ii \R \subseteq \{0\}$ already for $\Re \sp{\mathcal{A} \hat x}{\hat x}_{\mathcal{X}} \leq - \kappa \abs{\mathfrak{R} x}^2$ where $\mathfrak{R} x = ((\H^1_2 x^1_2)(0), (\H^1 x^1)'(0))$.)
  The only candidate for an eigenfunction in this case satisfies
   \[
    (\H_2 x_2)(\z) = 0,
     \quad
     (\H_1 x_1)(\z) = \z (\H_1 x_1)'(0)
   \]
  i.e.\ from $(\H^j_1 x^j_1)(1) = (\H^{j+1}_1 x^{j+1}_1)(0)$ it follows that
   \[
    (\H^j_1 x^j_1)(\z)
     = (j - 1 + \z) (\H^1_1 x^1_1)'(0)
     = (j - 1 + \z) c,
     \quad \z \in [0,1], \, j = 1, \ldots, m
   \]
  for some $c \in \K$ and the eigenspace $\ker (\mathcal{A})$ is one-dimensional.
  This corresponds to the dynamical solution $\omega_t(t,\z) = (j - 1 + \z) c + \omega_t(0,\z)$, clearly a solution which already for moderately large $t > 0$ does not satisfy the underlying modelling assumptions for a linear Euler-Bernoulli beam model. In particular, the assumption that $\abs{\omega(t,\z) - \omega^{\mathrm{ref}}(\z)} \ll 1$ for some reference configuration $\omega^{\mathrm{ref}}$ will be violated then.
  A physical interpretation of this eigenstate would be a beam which is rotating in the transversal flat.
  Such phenomena are, as already stated, not covered by the linear beam model.
  However, after restricting the initial data to $\ker(\mathcal{A})^\perp$, i.e.\
   \begin{align*}
    \sum_{j=1}^m \int_0^1 (j-1+\z) (\H^j_1 x^j_1)(0,\z) \dd \z
     &= \sum_{j=1}^m \int_0^1 (j-1+\z) \tilde \H_1 ((1-\z) l^j + \z l^{j+1}) \tilde x_1((1-\z)l^j + \z l^{j+1})
     \\
     &= \sum_{j=1}^m \int_{l^{j-1}}^{l^j} \frac{j - 1 + \z}{\z_j - \z_{j-1}} \omega_t(0,\z)
     = 0,
   \end{align*}
  by linearity, compactness of the resolvent ($\sigma_p(\mathcal{A}) \cap B_1(0)$ is discrete!) and Theorem \ref{thm:exp_stability} also in the second and third case, the solution tends uniformly exponentially to zero, even for the choice $K_0 = \operatorname{diag}\, (k^0_{11}, 0)$ for some $k^0_{11} > 0$.
 \end{rmrk}
 
\section{Conclusion}
\label{conclusion}

In this paper we presented a proof via the resolvent method for the uniform stabilisation of a chain of serially connected inhomogeneous Euler-Bernoulli beams with damping at one end.
We considered several possible interconnection conditions and pairs of dissipative / conservative boundary conditions at the ends of the chain which enforce uniform exponential energy decay for the beam system.
We thereby not only generalised the results in \cite{ChenEtAl_1987} to the case of non-uniform beams (which in this generality seems not to be possible by their method), but identified several other possible combinations of dissipative-conservative pairs boundary conditions at the left and right end of the chain leading to exponential energy decay as well.
Moreover, we showed that instead of static boundary or feedback interconnections, dynamic feedback interconnections with finite dimensional control systems can be used as well to achieve well-posedness and stability results.

\section*{Acknowledgment}
I am much obliged to Birgit Jacob who introduced me to the topic of infinite-dimensional port-Hamiltonian systems and encouraged me to carry on research in this area. Moreover, I would like to thank the anonymous referees for their careful reading and their valuable remarks, leading to significant improvement of the manuscript.

\end{document}